\documentclass{article}
\usepackage{epsfig}
\usepackage{latexsym}
\usepackage{amsmath}
\usepackage{amssymb}
\usepackage[english]{babel}

\newtheorem{thm}{Theorem}[subsection]
\newtheorem{cor}[thm]{Corollary}
\newtheorem{lem}[thm]{Lemma}
\newtheorem{prop}[thm]{Proposition}

\newtheorem{expl}[thm]{Example}

\newcommand{\M}{\mathcal{M}}

\newcommand{\Real}{\mathbb{R}}

\DeclareMathOperator*{\bigboxplus}{\text{\Huge $_{\boxplus}$}}

\DeclareMathOperator*{\medotimes}{\text{\huge $_{\otimes}$}}
\DeclareMathOperator*{\bigsmileplus}{\text{\Large ${\stackrel{_{+}}{\smile}}$}}
\DeclareMathOperator*{\bigsmileminus}{\text{\Large ${\stackrel{-}{\smile}}$}}

\DeclareMathOperator*{\Card}{\mathrm{Card}}
\DeclareMathOperator*{\sgn}{\mathrm{sgn}}

\DeclareMathOperator*{\cl}{\mathrm{cl}}
\DeclareMathOperator*{\inter}{\mathrm{int}}

\numberwithin{equation}{section}

\newfont{\sfl}{cmssi12}

\begin{document}

\title{ Determinants and Limit Systems
in some Idempotent  and Non-Associative  Algebraic
Structure}

 \author{\thanks{University of Perpignan, LAMPS, 52 avenue Villeneuve,
66000 Perpignan, France. email{briec@univ-perp.fr} }Walter Briec}
\date{September, 2020}
 \maketitle

\begin{abstract}
 This paper considers  an idempotent and
 symmetrical algebraic structure as well as some closely related concept. A special notion of determinant is introduced  and a Cramer formula is derived for a class of limit systems derived from the Hadamard matrix product and we give the algebraic form of a sequence of hyperplanes passing through a finite number of points.   Thereby, some standard results arising for Max-Times systems   with nonnegative entries appear as  a special case. The case of two sided systems is also analyzed. In addition, a notion of eigenvalue in limit is  considered. It is shown that one can construct a special  semi-continuous regularized polynomial to find the eigenvalues of a matrix with nonnegative entries.
\end{abstract}

{\bf AMS:} 06D50, 06F25\\


\section{Introduction}\label{SECMXVSP}

Exotic or tropical semirings such as the Max-Plus semiring, have been developed since the late fifties. They have many applications to various fields:
performance evaluation of manufacturing systems;  graph
theory and Markov decision processes; Hamilton-Jacobi theory.
However, it is well known that there is no nontrivial algebraic structures satisfying both idempotence, symmetry and having a neutral element. Despite this, there exist methods for symmetrizing  an idempotent semiring  imitating the familiar construction of $\mathbb Z$
from $\mathbb N$, for an arbitrary semiring. Symmetrization of  idempotent Semirings plays a crucial role to
develop an approach in term of determinant in Max-Plus Algebra.  Gaubert \cite{gau92}
introduced a balance relation to preserve transitivity. Familiar identities valid in rings
admit analogues, replacing equalities by balances. The balance relation
yields to relations similar to those arising for ordinary
determinant making a lexical change. This symmetrization was
invented independently by G. Heged\"us \cite{h85} and M. Plus
\cite{MP}. It follows that solving linear equations in the Max-Plus semi-ring requires to
solve systems of linear balances. Results concerning Cramer
solutions can be found in   \cite{BCOQ92}.

In this paper we have taken a different point of view. We consider an idempotent algebraic structure having the symmetry property and $0$ as a neutral element. The price to pay is that associativity no longer holds true. More precisely, we focus on a Max-Times algebraic structure which is derived as a limit case of the generalized power-mean involving an homeomorphic transformation of the real field. The binary operation involved by this algebraic structure was mentioned in \cite{Gau98} as an exercice.\footnote{Exercice 41, p. 25.}   Though it is not associative it admits an $n$-ary extension and satisfies some interesting properties. In particular, one  can construct a scalar product which will play an important role in the paper. It has been shown  in \cite{b15,b17} that such an algebraic structure is useful to extend a Max-Times idempotent convex structure from $\Real_+^n$ to the whole Euclidean vector vector space.  The problem arising with such a cancelative algebraic structure  is that it involves a natural $n$-ary operation  that is not continuous nor associative. Therefore, to circumvent this difficulty and establish separation properties of convex sets,  a special class of semi-continuous (upper and lower) regularized inner products  was considered in \cite{b17}.

The paper focusses on the asymptotic Cramer solutions of a special  sequence of  generalized power-linear systems.   These systems are constructed from an homeomorphic transformation  of the usual matrix product  involving the Hadamard power for vectors and matrices.\footnote{ A similar approach was considered in \cite{abg98} modulo a logarithmic change in the variables related to the Max-Plus algebraic structure.}  The formula of the determinant and Cramer's rule are then derived with respect to the non-associative algebraic structure considered in  \cite{b15}. Along this line, we give the algebraic form of a sequence of hyperplanes passing through a finite number of points.  More importantly, a general class of limit systems is defined over $\Real^n$.  These limit systems involve  several inequations that are derived from the semi-continuous (upper and lower) regularization of the non-associative  inner product. They includes as a special case   all the Max-Times systems defined from a matrix  with positive entries. The Kaykobad's conditions established in \cite{kayko} can then be applied to warrant the asymptotic existence of a positive solution.   This algebraic structure does not require any balance relation and  one can give an explicit form to some solutions of a two-sided Max-Times system. In addition, it is shown that one can construct a special  polynomial to find the eigenvalues of a matrix with nonnegative entries. To do that the limit of the Perron-Frobenius eigenvalue is considered. A parallel viewpoint was adopted in \cite{abg98} in a Max-Plus context.

The paper unfolds as follows. We lay down the groundwork is section 2.  In section 3, a suitable notion of determinant is   defined with respect to this non-associative algebraic structure. Section 4 considers a class of semi-continuous regularized operators. Hence an explicit algebraic form of the limit of a sequence of generalized hyperplanes is provided. In section 5 a class of limit systems of equations   is analyzed for which an explicit  Cramer formula is established including the case of Max-Times systems with nonnegative entries. In addition, we provide a solution for a class of two-sided systems and we compare the balance relations and the non-associative algebraic structure used in the paper. Finally, a notion of eigenvalues in limit is analyzed and connected to the  algebraic structure proposed in the paper.

\section{Preliminary Properties}

\subsection{An Idempotent and Non-Associative Algebraic Structure}
For all $p\in \mathbb{N}$,
let us consider a bijection $\varphi_p :\mathbb{R} \longrightarrow
\mathbb{R}$ defined by:
\begin{equation}\varphi_p:x \longrightarrow x^{2p+1}\end{equation}
 and $\phi_p (x_1,... ,x_n)=(\varphi_p(x_1),... ,\varphi_p(x_n))$; this is closely related to the approach proposed by
Ben-Tal   \cite{ben} and Avriel   \cite{avr1}. One can
induce a field structure on $\mathbb{R}$ for which $\varphi_p$
becomes a field isomorphism. Given this change of notation via
$\varphi_p$ and $\phi_p$ we can define a $\Real$-vector space
structure on $\Real^n$ by: $\lambda \stackrel{\varphi_p}{.}
x=\phi_p^{-1}(\varphi_p(\lambda).\phi_p(x))=\lambda.x$ and $x
\stackrel{\varphi_p}{+} y=\phi_p^{-1}(\phi_p(x)+\phi_p(y))$;
we call these two operations the indexed scalar product and the
indexed sum (indexed by $\varphi_p$).

The $\varphi_p$-sum  denoted $\stackrel{\varphi_p}{\sum}$    of
$(x_1,...,x_m)\in \mathbb{R}^{n\times m}$ is defined by\footnote{For all positive natural numbers  $n$, $[n]=\{1,...,n\}$.}
 \begin{equation}\stackrel{\varphi_p}{\sum_{i\in [m]}}x_j=
\phi_p^{-1}\Big(\sum_{j\in [m]}\phi_p (x_j)\Big).\end{equation} For simplicity, throughout the paper we
denote for all $x,y\in \Real^n$:
\begin{equation}x\stackrel{p}{+}y=x\stackrel{\varphi_p}{+}y.\label{phip}\end{equation}
Recall that  Kuratowski-Painlev\'e lower limit of the sequence
of sets $\{A_n\}_{n\in \mathbb N}$, denoted $Li_{n\to\infty}A_n$, is the set of
points $x$ for which there exists a sequence $\{x^{(n)}\}_{n\in \mathbb N}$ of points
such that $x^{(n)}\in A_n$ for all $n$ and $x = \lim_{n\to\infty}x^{(n)}$. The Kuratowski-Painlev\'e upper limit of the sequence
of sets $\{A_n\}_{n\in \mathbb N}$, denoted $Ls_{n\to\infty}A_n$, is the set of
points $x$ for which there exists a subsequence $\{x_{n_k}\}_{k\in \mathbb N}$ of points
such that $x^{(n_k)}\in A_{n_k}$ for all $k$ and $x = \lim_{k\to\infty}x^{(n_k)}$.
A sequence $\{A_n\}_{n\in \mathbb{N}}$ of subsets of $\Real^n$ is
said to converge, in the Kuratowski-Painlev\'e sense, to a set $A$
if $Ls_{n\to\infty}A_n = A = Li_{n\to\infty}A_n$, in which case we
write $A = Lim_{n\to\infty}A_n$.

\subsection{ A Limit Algebraic Structure}\label{Rec}

 In \cite{b15} it was shown that for all $x,y\in \Real$ we have:
\begin{equation*}\label{base}\lim_{p\longrightarrow
+\infty}x\stackrel{p}{+}y=
\left\{\begin{matrix}x\ &\hbox{ if } &|x|&>&|y|\\
\frac{1}{2}(x+y)&\hbox{ if }&|x|&=&|y|\\
y& \hbox{ if }& |x|&<&|y|.\end{matrix}\right.\end{equation*} Along
this line one can introduce the binary operation $\boxplus$
defined for all $x,y\in \Real$ by:
\begin{equation}
x\boxplus y=\lim_{p\longrightarrow +\infty}x\stackrel{p}{+}y.
\end{equation}
Though the operation $\boxplus$ does not satisfy associativity, it
can be extended by constructing a non-associative algebraic
structure which returns to a given $n$-tuple a real value. For all
$x\in \mathbb R^n$ and all subsets $I$ of $[n]$, let us  consider
the map $\xi_I[x]:\Real \longrightarrow \mathbb Z$ defined for all
$\alpha \in \Real$ by
\begin{equation}\label{defxi}\xi_I[x](\alpha)= \Card
\{i\in I: x_i=\alpha\}- \Card \{i\in I: x_i=-\alpha\}.\end{equation}
This map measures the symmetry of the occurrences of a given value
$\alpha $ in the components of a vector $x$.

For all $x\in \mathbb R^n$ let $\mathcal J_I(x)$ be a subset of $I$
defined by
\begin{equation}\mathcal J_I(x)=\Big\{j\in I: \xi_I[x](x_j)\not=0\Big\}=I\backslash \big \{i\in I: \xi_I[x]=0\big\}.\end{equation}
$\mathcal J_{I} (x)$ is called {\bf the residual  index set } of
$x$. It is obtained by dropping from $I$ all the $i$'s such that
$\Card \{j\in I: x_j=x_i\}= \Card \{j\in I: x_j=-x_i\}$.

For all positive natural numbers $n$ and
for all subsets $I$ of $[n]$, let $\digamma_I: \Real^n
\longrightarrow \Real $ be the map defined for all $x\in \Real^n$ by

\begin{equation}\digamma_{ I}(x)=\left\{\begin{matrix}\max_{i\in \mathcal
J_I (x)}x_{i} &\text{ if }&\xi_I[x]\big (\max_{i\in
\mathcal J_I(x)}|x_i|\big)>0\\
\min_{i\in \mathcal J_I(x)}x_i &\hbox{ if }&\xi_I[x]\big(\max_{i\in \mathcal J_I(x)}|x_i|\big)<0\\
 0 &\text{ if }&\xi_I[x]\big (\max_{i\in
\mathcal J_I(x)}|x_i|\big)=0.
\end{matrix}\right.\end{equation}
where $\xi_I[x]$ is the map defined in \eqref{defxi} and $\mathcal
J_I(x)$ is the residual index set of $x$. The operation that takes
an $n$-tuple $(x_1,....,x_n)$ of $\Real^n$ and returns a single
real element $\digamma_I(x_1,...,x_n)$ is called a $n$-ary
extension of the binary operation $\boxplus$ for all natural
numbers $n\geq 1$ and all $x\in \Real^{n}$, if $I$ is a nonempty
subset of $[n]$. Then, for all $n$-tuple $x=(x_1,...,x_n)$, one
can define the operation:
\begin{equation}\label{defrecOp}\bigboxplus_{i\in
I}x_i= \lim_{p\longrightarrow
\infty}\stackrel{\varphi_p}{\sum_{i\in I}}x_i=\digamma_{I}(x).
\end{equation}

Clearly, this operation encompasses as a special case the binary
operation defined in equation \eqref{base} and for all
$(x_1,x_2)\in \Real^2$:
\begin{equation*}\bigboxplus_{i\in \{1,2\}}x_i=x_1\boxplus
x_2.\end{equation*} For example,  if $x=(-3,-2,3,3,1,-3)$,   we
have $F_{[6]}(-3,-2,3,3,1,-3)=F_{[2]}(-2,1)=-2=\bigboxplus_{i\in
[6]}x_i$. There are some basic properties that can be inherited
from the above algebraic structure. We briefly summarize some
basic properties: $(i)$ If all the elements of the family
$\{x_i\}_{i\in I}$ are mutually non symmetrical, then:
$\bigboxplus_{i\in I}x_i=\arg\max_{\lambda }\big\{|\lambda|:
\lambda \in \{x_i\}_{i\in I}\big\}$; $(ii)$ For all $\alpha\in
\Real$, one has: $ \alpha \Big(\bigboxplus_{i\in I}x_i\Big)=
\bigboxplus_{i\in I}(\alpha x_i)$; $(iii)$ Suppose that
$x\in\epsilon \Real_+^n$ where $\epsilon$ is $+1$ or $-1$. Then
$\bigboxplus_{i\in I}x_i=\epsilon \max_{i\in I} \{\epsilon x_i\}$;
$(iv)$ We have $|\bigboxplus_{i\in I} x_i|\leq \bigboxplus_{i\in
I}|x_i|$; $(v)$ For all $x\in \Real^n$:
$$\Big[x_i \boxplus \big(\bigboxplus_{j\in I\backslash
\{i\}}x_j\big)\Big]\in \Big \{0, \bigboxplus_{j\in I}x_j\}\quad
\text{ and }\quad \bigboxplus_{i\in I}x_i=\bigboxplus_{i\in
I}\Big[x_i \boxplus \big(\bigboxplus_{j\in I\backslash
\{i\}}x_j\big)\Big].$$

The algebraic structure $(\Real,\boxplus,\cdot)$  can be extended
to $\Real^n$. Suppose that $x,y\in \Real^n$, and let us denote $
x\boxplus y=(x_1\boxplus y_1,...,x_n\boxplus y_n).$  Moreover,
let us consider $m$ vectors $x_1,...,x_m\in \Real^n$, and define
\begin{align}
\bigboxplus_{j\in [m]}x_j&=\Big(\bigboxplus_{j\in [m]}
x_{j,1},...,\bigboxplus_{j\in [m]} x_{j,n}\Big).
\end{align}

The $n$-ary operation $(x_1,...,x_n)\rightarrow \bigboxplus_{i\in [n]}x_i$ is not associative.
To simplify the notations of the paper, for  all $z\in \{z_{i_1,...,i_m}: i_k\in I_k,k\in [m]\}$, where $I_1,...,I_m$ are $m$ index subsets of $\mathbb N$, we use the notation:
\begin{equation}
\bigboxplus_{ \substack{i_k\in I_k\\k\in [m]}}z_{i_1,...,i_m}=\bigboxplus_{\substack{(i_1,...,i_m)\in \prod_{k\in [m]} I_k}}z_{i_1,...,i_m}.
\end{equation}
Notice that for all $x\in \Real^n$ and all $y\in \Real^m$:
\begin{equation}\label{product}
\Big(\bigboxplus_{i\in [n]}x_i\Big)\Big(\bigboxplus_{j\in [m]}y_j\Big)=\bigboxplus_{ \substack{i\in [n]\\j\in [m]}}x_iy_j.
\end{equation}
This relation immediately comes from the fact that for all natural numbers $p$, we have:
\begin{equation}
\Big(\stackrel{\varphi_p}{\sum_{i\in [n]}}x_i\Big)\Big(\stackrel{\varphi_p}{\sum_{j\in [m]}}y_j\Big)=\stackrel{\varphi_p}{\sum_{i\in [n]}}\stackrel{\varphi_p}{\sum_{j\in [m]}}x_i y_j=\stackrel{\varphi_p}{\sum_{ \substack{i\in [n]\\j\in [m]}}}x_iy_j.
\end{equation}
Taking the limit on both sides yields equation \eqref{product}.  In the remainder, we will adopt the following notational convention. For all $x\in \Real^n$:
\begin{equation}
\bigboxplus_{i\in [n]}x_i=x_1\boxplus\cdots\boxplus x_n=\digamma_{[n]}(x).
\end{equation}

\subsection{ Scalar Product}

\noindent This section presents the algebraic properties induced
by the isomorphism of scalar field $\varphi_p$  on the scalar product. Most of the results have been pointed
in details by Avriel \cite{avr1} and Ben Tal \cite{ben}. A norm
$\|\cdot\|$ yields another norm induced by the algebraic
operations $\stackrel{p}{+}$ and $\cdot$. The map $\|\cdot
\|_{\varphi_p}:\Real^n\longrightarrow \Real$ defined by
$\|x\|_{\varphi_p}=\varphi_p^{-1}\left(\|\phi_p(x)\|\right)$ is a
norm over $\Real^{n}$ endowed with the operations
$\stackrel{p}{+}$ and $\cdot$. Since $\varphi_p$ is continuous
over $\mathbb{R}$, the topological structure is the same. Along
this line it is natural to define a scalar product. If $\langle
\cdot ,\cdot \rangle$ is an inner product over $\Real^{n}$, then
there exists a symmetric bilinear form
 $\langle \cdot,\cdot \rangle_{\varphi_p}:
 \Real^n\times \Real^n\longrightarrow \Real$ defined by:
\begin{equation}\langle x,y\rangle_{\varphi_p}=\varphi_p^{-1}\big(\langle \phi_p^{}(x),
\phi_p^{}(y)\rangle\big )=\big(\sum_{i\in [n]}x_i^{2p+1} y_i^{2p+1}\big)^{\frac{1}{2p+1}}.\end{equation}
Now, let us denote $\big [\langle y,\cdot\rangle_{\varphi_p}\leq
\lambda\big ]= \left\{x\in \Real^n: \langle
y,x\rangle_{\varphi_p}\leq \lambda\right\}$ and let  $\langle \cdot,\cdot \rangle_p$ stands for this scalar product.

In the following we introduce the operation $ \langle \cdot,\cdot\rangle_\infty :\Real^n\times \Real^n\longrightarrow
 \Real$ defined for all $x,y\in \Real^n$ by $\langle x,y\rangle_\infty =\bigboxplus_{i\in
 [n]}x_iy_i$. Let $\|\cdot\|_\infty$ be the Tchebychev
 norm defined by $\|x\|_\infty=\max_{i\in [n]}|x_i|$. It is established in \cite{b15} that
for all $x,y\in \Real^n$, we have:
 $(i)$  $\sqrt{\langle x,x\rangle_\infty} =\|x\|_\infty$;
  $(ii)$  $|\langle x,y\rangle_\infty| \leq \|x\|_\infty \|y\|_\infty$;
   $(iii)$ For all $\alpha\in \Real$, $\alpha \langle x,y\rangle_\infty= \langle \alpha x,y\rangle_\infty=\langle  x,\alpha y\rangle_\infty$.
By definition, we have for all $x,y\in \Real^n$:

\begin{equation}\langle x,y\rangle_\infty=\lim_{p\longrightarrow \infty}\langle x,y\rangle_p\end{equation}

\section{ Limit of Linear Operators and Determinant}

This section is devoted to study the matrix representation of a
linear operator defined on the scalar field $(\Real,\stackrel{p}{+},\cdot)$. Along this line some limit properties are derived
 to establish several results in closed algebraic form
when $p\longrightarrow \infty$.

\subsection{$\varphi_p$-linear Endomorphisms}
  Let $\mathcal  L(\Real^n, \Real^m)$ denotes
 the set of all the linear endomorphisms defined from  $\Real^n$ to $\Real^m$.
 Let $\mathcal  L(\Real^n, \Real^n)$ is then
 the set of all the linear endomorphisms defined over  $\Real^n$.
In the following, we say that a map $f: \Real^n\longrightarrow
\Real^n$ is $\varphi_p$-linear if for all $\lambda\in \Real$,
$f(\lambda x\stackrel{p}{+} y)=\lambda f( x)\stackrel{p}{+}f(y) $.
Moreover, for all natural numbers $p$, let  $\mathcal  L
^{(p)}(\Real^n, \Real^n)$ denotes
 the set of all  the $\varphi_p$-linear endomorphisms.

 Let $\M_{n}(\Real)$ denotes the set of all the
$n\times n$ matrices defined over  $\Real$. Let $\Phi_p:\M_{n}(\Real)\longrightarrow \M_{n}(\Real)$ be the map defined
for any matrix $A= \left(a_{i,j}\right)_{\substack
{i=1...n\\j=1..n}} \in \M_{n}(\Real)$ as:
\begin{align}
\Phi_p(A)&=\big(\varphi_p(a_{i,j})\big)_{\substack
{i=1...n\\j=1...n}}= \left({a_{i,j}}^{2p+1}\right)_{\substack
{i=1...n\\j=1...n}}.\end{align} Its reciprocal is the map
$\Phi_p^{-1}:\M_{n}(\Real)\longrightarrow \M_{n}(\Real)$ defined
by:\begin{align}
\Phi_p^{-1}(A)&=\left(\varphi_p^{-1}(a_{i,j})\right)_{\substack
{i=1...n\\j=1...n}}=\left({a_{i,j}}^{\frac{1}{2p+1}}\right)_{\substack
{i=1...n\\j=1...n}}.
\end{align}
$\Phi_p$ is a natural extension of the map $\phi_p$ from $\Real^n$
to $\mathcal M_n(\Real)$. $\Phi_p(A)$ is the $2p+1$ Hadamard power
of matrix $A$. In the following we introduce the matrix
product:
\begin{equation}
A \stackrel{p}{.}x=\sum_{j\in [n]}^{\varphi_p} x_{j}{.}a^{j},
\end{equation}
where $a^j$ stands for the $j$-th column of $A$. It is straightforward to
show that this formulation is equivalent to the following:
\begin{equation}
A \stackrel{p}{.}x=\phi_p^{-1}\big(\Phi_p(A).\phi_p(x)\big).
\end{equation}
 Another equivalent formulation involves the inner product $\langle
\cdot,\cdot\rangle_p$:

\begin{equation}
A \stackrel{p}{.}x=\sum_{i\in [n]}^{\varphi_p} \langle
a_{i},x\rangle_{p}e_{i},
\end{equation}
where $a_i$ is the $i$-th line of matrix $A$ and $\{e_i\}_{i\in
[n]}$ is the canonical basis of $\Real^n$.

It is easy to see that the map $x\mapsto A\stackrel{p}{\cdot}x$ is
$\varphi_p$-linear. Conversely,  if $g$ is a $\varphi_p$-linear
map then it can be represented by a matrix $A$ such that $g(x)=
A\stackrel{p}{\cdot}x$ for all $x\in \Real^n$. If $A,B\in \mathcal
M_n(\Real)$, the product $A\stackrel{p}{\cdot}B$ is the matrix representation of the map:
\begin{equation} x\mapsto B\stackrel{p}{\cdot}A\stackrel{p}{\cdot}x=\phi_p^{-1}\big(\Phi_p(B)\Phi_p(A)\phi_p(x)\big).\end{equation}
Notice that the identity matrix $I$ is invariant with respect to
$\Phi_p$.

Let $f:\Real^n\longrightarrow \Real^n$ be a linear endomorphism and let $A$ be its matrix representation in the canonical
basis. The map $T^{(p)}:\mathcal L(\Real^n, \Real^n)\longrightarrow
\mathcal L ^{(p)}(\Real^n,\Real^n)$ defined for all $x\in \Real^n$ by:
$$T^{(p)}(f)(x)=f^{(p)}(x):=\phi_p^{-1}\big(\Phi_p(A)\phi_p(x)\big)$$
is called {\bf the $\varphi_p$-linear transformation of $f$}.

A $\varphi_p$-linear endomorphism $g$ is invertible if and only if
$\Phi_p(A)$ is invertible.   For any $n\times n$ matrix $A$,
let $|A|$ denotes its determinant. Let us introduce the following
definition of a {\bf $\varphi_p$-determinant}
\begin{equation}|A|_p=\varphi_p^{-1}|\Phi_p(A)|.\end{equation}
Let $S_n$ be the set of all the permutations defined on $[n]$. The Leibnitz formula yields
\begin{equation}|A|_p=\Big(\sum_{\sigma\in S_n}\sgn(\sigma) \prod_{i\in [n]}a_{i,\sigma(i)}^{2p+1}\Big)^{\frac{1}{2p+1}}.\end{equation}

In the remainder, if $A$ is the matrix of a linear endomorphism $f$, then we define the $\varphi_p$-determinant of $f$ as $|f|_p=|A|_p$. If $f^{(p)}$ is invertible, then we have the equivalences:
 \begin{equation}y=f^{(p)}(x)\iff y=\phi_p^{-1}\big(\Phi_p(A)\phi_p(x)\big)\iff \phi_p^{-1}\big({\Phi_p(A)}^{-1}\phi_p(y)\big)=x.\end{equation}Along this line, the $\varphi_p$-cofactor matrix
$A^{\star,p}$ is defined as:
\begin{equation}
A^{\star,p}=(a^{\star,p}_{i,j})_{\substack{i\in [n]\\j\in
[n]}}=\Big((-1)^{i+j}\big|A_{i,j}\big|_p\Big)_{\substack{i\in
[n]\\j\in [n]}},
\end{equation}
 where $A_{i,j}$ is obtained from matrix
 $A$ by dropping line $i$ and column $j$. The
 $\phi_p$-inverse matrix of a $\varphi_p$-invertible matrix $A$ (such that $|A|_p\not=0$) is then defined as:
 \begin{equation}A^{-1,p}=\frac{1}{|A|_p} { ^t\!A}^{\star,p}.\end{equation}
Suppose that $f$ is a linear endomorphism having a matrix
representation $A$ in the canonical basis and let $b\in \Real^n$.
Given a system of $\varphi_p$-linear equations of the form:
\begin{equation}f^{(p)}(x)=b\iff A\stackrel{p}{\cdot}x=b,\end{equation}
if $|A|_p\not=0$, then the solution is $x^\star=A^{-1,p} =b$.

\subsection{Limit Properties}

\begin{prop}  Let $f:\Real^n\longrightarrow \Real^n$ be a linear
endomorphism having a matrix representation $A$. For all $p\in
\mathbb N$, let $f^{(p)}$ be its $\varphi_p$-linear
transformation. Then:
$$\lim_{p\longrightarrow \infty}f^{(p)} (x) =\bigboxplus_{j\in [n]}x_ja^j
=\bigboxplus_{i\in [n]}\langle a_i, x\rangle_\infty e_i.$$
 \end{prop}
 {\bf Proof:}  For all $i \in [n]$, we have from \cite{b15}:
 $$\lim_{p\longrightarrow \infty}\stackrel{\varphi_p}{\sum_{j\in [n]}}{x_j}a_ {i,j}
 =\lim_{p\longrightarrow \infty}\big({\sum_{j\in [n]}}{x_j}^{2p+1}a_ {i,j}^{2p+1}\big)^{\frac{1}{2p+1}}=
 \bigboxplus_{j\in [n]}x_ja_{i,j} =\langle a_i, x\rangle_\infty.$$
 Therefore $$\lim_{p\longrightarrow \infty}\stackrel{\varphi_p}{\sum_{j\in [n]}}{x_j}a^j=\bigboxplus_{j\in [n]}x_ja^j.$$
 The last equality immediately follows. $\Box$

For any squared matrix $A$, $|A|_\infty$ is called {\bf the determinant in limit of $A$}. For any linear endomorphism $f$ whose the matrix is $A$, the determinant in limit of $f$ is defined as $|f|_\infty=|A|_\infty. $
\begin{prop} \label{limdet}  For all $A\in \mathcal M_n(\Real)$, we have:
$$\lim_{p\longrightarrow \infty}|A|_p: =|A|_\infty=\bigboxplus_{\sigma\in S_n} \big(\sgn(\sigma)
 \prod_{i\in [n]}a_{i,\sigma(i)}\big).$$
\end{prop}
{\bf Proof:} From \cite{b15} we have:
$$\lim_{p\longrightarrow \infty}\Big(\sum_{\sigma\in S_n}\sgn(\sigma)
\prod_{i\in [n]}a_{i,\sigma(i)}^{2p+1}\Big)^{\frac{1}{2p+1}}
=\bigboxplus_{\sigma\in S_n} \big(\sgn(\sigma)
 \prod_{i\in [n]}a_{i,\sigma(i)}\big).\quad \Box$$

 \begin{prop} \label{limsol1} Let $f:\Real^n\longrightarrow \Real^n$ be a linear
endomorphism having a matrix representation $A$. For all $p\in
\mathbb N$, let $f^{(p)}$ be its $\varphi_p$-linear
transformation.   If $|A|_\infty\not=0$, then there is some
$p_0\in \mathbb N$ such that  for all $p\geq p_0$,  $f^{(p)}$ is
$\varphi_p$-invertible and for all $b\in \Real^n$, there
exists a solution $x^{(p)}$ to the system
$A\stackrel{p}{\cdot}x=b$  with:
$$x^{(p)}_i=\frac{|A^{(i)}|_p}{|A|_p}=\frac{\big |\Phi_p(A^{(i)})\big |^{\frac{1}{2p+1}}}{\big |\Phi_p(A)\big |\frac{1}{2p+1}},$$
where $A^{(i)}$ is obtained from $A$ by dropping column $i$ and
replacing it with $b$. Moreover, we have:
  $$ \lim_{p\longrightarrow \infty}x^{(p)}=x^\star,$$
  with for all $i\in[n]$
  $$x^\star_i=\frac{|A^{(i)}|_\infty}{|A|_\infty}.$$
  \end{prop}
  {\bf Proof:} Since $|A|_\infty\not=0$ and $\lim_{p\longrightarrow
  \infty}|A|_p=|A|_\infty$, there is some $p_0$ such that for all
  $p\geq p_0$, $|A|_p\not=0$, which implies that $A$ is
  $\varphi_p$-invertible. In such a case, there exists an
  uniqueness solution to the system $A\stackrel{p}{\cdot}x=b$,
 that is $x^{(p)}=A^{-1,p}b$. Moreover, we have:
$$A\stackrel{p}{.}x=b\iff \phi_p^{-1}\big(\Phi_p(A)\phi_p(x)\big)=b\iff \Phi_p(A)\phi_p(x)=\phi_p(b).  $$
Since $f^{(p)}$ is $\varphi_p$-invertible, it follows that
$\Phi_p(A)$ is invertible. Set $u=\phi_p(x)$. The system
$\Phi_p(A)u=\phi_p(b)$ has a solution for all $p\geq p_0$.
Applying the Cramer's rule the solution  is the vector $u^{p}$
satisfying the relation:
$$u^{(p)}=\frac{\big|[\Phi_p(A)]^{(i)}\big | }{\big |\Phi_p(A)\big | }=\frac{\big |\Phi_p(A^{(i)})\big | }{\big |\Phi_p(A)\big | }.$$
Setting $x^{(p)}=\phi_p^{-1}(u^{(p)})$, we obtain the result. From
Proposition \ref{limdet}, $\lim_{\longrightarrow
\infty}|A^{(i)}|_p=|A^{(i)}|_\infty$ and $\lim_{\longrightarrow
\infty}|A|_p=|A|_\infty$, which ends the proof. $\Box$\\

The next properties are useful. We first establish the following Lemma.

\begin{lem} \label{induc}Suppose that there is some $x=(x_1,...,x_n)\in \Real^n$ such that $\bigboxplus_{i\in [n]}x_i=0$. Then for all $p\in \mathbb N$, $\stackrel{\varphi_p}{\sum_{i\in [n]}}x_i=0$. Moreover, for all matrices $A\in \mathcal M_n(\Real)$, if $|A|_\infty=0$ then $|A|_p=0$ for all $p\in \mathbb N$.
\end{lem}
{\bf Proof:} Let $\Lambda[x]=\{\alpha \in \Real_+ :|x_i|=\alpha, i\in [n]\}$. Since $\bigboxplus_{i\in [n]}x_i=0$, we have  for all $\alpha$, $\xi[x](\alpha)
=\Card\{i: x_i=\alpha\}-\Card\{i: x_i=-\alpha\}=0$. Hence, for all $\alpha\in \Lambda[x]$, $\stackrel{\varphi_p}{\sum_{|x_i|=\alpha}}x_i=0$. Thus
$$\stackrel{\varphi_p}{\sum_{i\in [n]}}x_i=\stackrel{\varphi_p}{\sum_{\alpha\in A[x]}}\;\stackrel{\varphi_p}{\sum_{|x_i|=\alpha}}x_i=0.$$
The second part of the statement is an immediate consequence of the Leibniz formula. $\Box$\\

For all $p\in \mathbb N\cup\{\infty\}$ and for all matrices $A\in \mathcal M_n(\Real)$, let us denote  $|a^1,...,a^n|_p=|A|_p$ where the $a^j$'s are the column vectors of $A$.
\begin{prop}For all $A\in \mathcal M_n(\Real^n)$, we have the following properties.\\
$(a)$ For all $\alpha\in \Real$,   $| a^{ 1 },...,\alpha a^{ j },...,a^{ n } |_\infty=\alpha |A|_\infty$;\\
$(b)$ For all permutations $\sigma $ of $S_n$,  $|a^{ \sigma (1) },...,a^{ \sigma(n) } |_\infty=\mathrm{sgn} (\sigma) |A|_\infty$;\\
$(c)$ If there exists $\alpha\in \Real^n\backslash\{0\}$ such that
 $\bigboxplus_{j\in [n]}\alpha_ja^{ j }=0$  then $|A|_\infty=0$.\\
$(d)$ If  $|A|_\infty=0$ then there exists a sequence $\{\alpha^{(p)}\}_{p\in \mathbb N}\subset \Real^n\backslash\{0\}$ such that $\stackrel{p}{\sum_{j\in [n]}}\alpha_j^{(p)}a^{ j }=0$ for all $p\in\mathbb N$.
\end{prop}
{\bf Proof:} $(a)$   Since $ |A |_p=\varphi_p^{-1}\big(|\Phi_p(A)|\big)$, we deduce that for all natural numbers $p$, $|a^{ 1 },...,\alpha a^{ j },...,a^{ n }|_p=\alpha |a^{ 1 },...,a^{ n }|_p$. Taking the limit yields the result. $(b)$ Similarly, for all  permutations $\sigma \in S_n$, $ |a^{ \sigma (1) },...,a^{ \sigma(n) }|_p=\mathrm{sgn} (\sigma) |A|_p$, which yields $(b)$, by taking the limit. $(c)$ If there exists $\alpha\in \Real^n\backslash\{0\}$ such that
 $\bigboxplus_{j\in [n]}\alpha_ja^{ j }=0$, from Lemma \ref{induc}, then we deduce that for all natural numbers $p$, $\stackrel{p}{\sum_{j\in [n]}}\alpha_j^{(p)}a^{ j }=0$. However, this implies that  $|A|_p=0$ for all $p$. Hence
 $|A|_\infty=\lim_{p\longrightarrow \infty}|A|_p=0$. $(d)$ If $|A|_\infty=0$, then
$\bigboxplus_{\sigma\in S_n} \big(\sgn(\sigma)
 \prod_{i\in [n]}a_{i,\sigma(i)}\big)=0$. Thus, from Lemma \ref{induc} $$\stackrel{\varphi_p}{\sum_{\sigma\in S_n}} \big(\sgn(\sigma)
 \prod_{i\in [n]}a_{i,\sigma(i)}\big)=0$$ for all natural numbers $p$. Hence, for all $p$, there is
 $\alpha^{(p)}\in  \Real^n\backslash\{0\}$ such that $\stackrel{p}{\sum_{j\in [n]}}\alpha_j^{(p)}a^{ j }=0$, which ends the proof. $\Box$\\

Determinants are intimately linked to the exterior product  product of vectors that is an algebraic construction used  to study areas, volumes, and their higher-dimensional analogues.   Paralleling the earlier definitions,
a map $f:\Real^n\longrightarrow \Real$ is called a {\bf $\varphi_p$-multilinear form} if it is $\varphi_p$-linear in each argument. A $\varphi_p$-multilinear form is alternating if for each permutation $\sigma\in S_n$ we have $f(x_1,...,x_n)=\mathrm{sgn}(\sigma)f(x_{\sigma(1)},...,x_{ \sigma(n)})$. For all natural numbers $r$ an alternating $\varphi_p$-linear  $r$-form is a map defined for all $x_1,x_2,...,x_r\in \Real^n$ as:
\begin{equation}
\big(f_1\stackrel{p}{\wedge} f_2\cdots \stackrel{p}{\wedge}f_r\big )(x_1,...,x_r)=\stackrel{\varphi_p}{\sum_{\sigma\in S_r}}\mathrm{sign}(\sigma)f_1^{(p)}(x_{\sigma(1)})\cdots f_n^{(p)}(x_{\sigma(r)}),
\end{equation}
where for any $i$, $f_i$ is a linear  form and, and $f_i^{(p)}$ is the corresponding $\varphi_p$-transformation.
$\stackrel{p}{\wedge}$ is called the $\varphi_p$-exterior product of the linear forms $f_1,...,f_r$.
Let $\{e_1^\star,...,e_n^\star\}$ be the canonical basis of the dual space $\mathcal L(\Real^n,\Real)$. Suppose that $r=n$ and let $f=\sum_{i\in [n]}f(e_i)e_i^*$ be the linear endomorphism constructed from $f_1,...,f_n.$

\begin{prop} Let us consider
$n$ linear forms $f_1,...,f_n$. Then for all $x_1,x_2,...,x_n\in \Real^n$,
we have
\begin{equation*}
 \big(f_1 \stackrel{p}{\wedge}f_2\cdots \stackrel{p}{\wedge}f_n \big)(x_1,...,x_n) =|f|_p\big(e_1^\star\stackrel{p}{\wedge} e^\star_2\cdots \stackrel{p}{\wedge}e^\star_n\big )(x_1,...,x_n).
\end{equation*}
 Moreover,  we have:
\begin{equation*}
\lim_{p\longrightarrow \infty}\big(f_1 \stackrel{p}{\wedge} f_2 \cdots \stackrel{p}{\wedge}f_n  \big )(x_1,...,x_n)=|f|_\infty\big(e_1^\star\stackrel{\infty}{\wedge} e^\star_2\cdots \stackrel{\infty}{\wedge}e^\star_n\big ),
\end{equation*}
where
\begin{align*}\big(e_1^\star\stackrel{\infty}{\wedge} e^\star_2\cdots \stackrel{\infty}{\wedge}e^\star_n\big )(x_1,...,x_n)=\lim_{p\longrightarrow \infty}( e_1^\star\stackrel{p}{\wedge} e_2^\star\cdots e_n^\star\stackrel{p}{\wedge} e_n^\star)( x_1 ,..., x_n )=|x_1,...,x_n|_\infty.\end{align*}
\end{prop}
{\bf Proof:} Suppose that for $i=1,...,n$ there is a vector $a_i\in \Real^n$
such that $f_i(x)=\langle a_i, x\rangle$. Then $f_i^{(p)}(x)= \varphi_p^{-1}\big(\langle \phi_p(a_i), \phi_p(x)\rangle\big)=\big(\sum_{i\in [n]}a_i^{2p+1}x_i^{2p+1}\big)^{\frac{1}{2p+1}}$. It follows that:
\begin{align*}
\big(f_1 \stackrel{p}{\wedge} f_2 \cdots \stackrel{p}{\wedge}f_n \big )(x_1,...,x_n)&= \Big({\sum_{\sigma\in S_n}}\mathrm{sign}(\sigma)\prod_{i\in [n]}\langle \phi_p(a_i), \phi_p(x_{\sigma(i)})\rangle  \Big)^{\frac{1}{2p+1}}\\&=\Big({\sum_{\sigma\in S_n}}\mathrm{sign}(\sigma)
\prod_{i\in [n]}a_i^{2p+1}x_{\sigma(i)}^{2p+1}\Big)^{\frac{1}{2p+1}} .
\end{align*}
For each $i$, let $g_i^{(p)}:\Real^n\longrightarrow \Real$ be the linear form defined  by $g_i^{(p)}(z)=\langle \phi_p(a_i), z\rangle$. It follows that:
\begin{equation*}
\big(f_1 \stackrel{p}{\wedge} f_2 \cdots \stackrel{p}{\wedge}f_n \big )(x_1,...,x_n)=
\varphi_p^{-1}\Big(\big( g_1^{(p)}\wedge g_2^{(p)}\cdots \wedge g_n^{(p)}\big)\big (\phi_p(x_1),...,\phi_p(x_n)\big)\Big) .
\end{equation*}
Let $\{e_1^\star,...,e_n^\star\}$ be the canonical basis of the dual space $\mathcal L(\Real^n,\Real)$. From the usual properties of an alternating  $n$-form we deduce that:
\begin{equation*}
\big( g_1^{(p)}\wedge g_2^{(p)}\cdots \wedge g_n^{(p)}\big)\big (\phi_p(x_1),...,\phi_p(x_n)\big)=|\Phi_p(A)|(e_1^\star\wedge e_2^\star\cdots \wedge e_n^\star)(\phi_p(x_1),...,\phi_p(x_n)),
\end{equation*}
where $\Phi_p(A)$ is the matrix whose line $i$ is the vector $\phi_p(a_i)$.
Since for all $i$ and all $x\in \Real^n$ we have $\langle e_i, x\rangle_p=x_i$, this canonical basis is also, independently  of $p$, the canonical basis of $\mathcal L_p(\Real^n,\Real)$. Since $|A|_p=\varphi_p^{-1}\big(|\Phi_p(A)|\big)$, it follows that:
\begin{align*}
\big(f_1^{(p)}\stackrel{p}{\wedge} f_2^{(p)}\cdots \stackrel{p}{\wedge}f_n^{(p)}\big )(x_1,...,x_n)&=| A |_p\Big(\sum_{i\in [n]}\langle e_i,\phi_p(x)\rangle \Big)^{\frac{1}{2p+1}} \\ &=| A |_p ( e_1^\star\stackrel{p}{\wedge} e_2^\star\cdots e_n^\star\stackrel{p}{\wedge} e_n^\star)( x_1 ,..., x_n ) .
\end{align*}
However, we have
$$( e_1^\star\stackrel{p}{\wedge} e_2^\star\cdots e_n^\star\stackrel{p}{\wedge} e_n^\star)( x_1 ,..., x_n )=| x_1,...,x_n |_p. $$
We then obtain the final result taking the limit. $\Box$\\

For all $(f_1,f_2,...,f_n)\in {\mathcal L(\Real^n, \Real)}^n$, let  $f_1\stackrel{\infty}{\wedge}f_2\cdots \stackrel{\infty}{\wedge}f_n$ denotes the pointwise limit of the sequence $\{f_1 \stackrel{p}{\wedge} f_2 \cdots \stackrel{p}{\wedge}f_n \}_{n\in \mathbb N}$. Namely
\begin{equation}
 f_1\stackrel{\infty}{\wedge}f_2\cdots \stackrel{\infty}{\wedge}f_n=\lim_{p\longrightarrow \infty}f_1 \stackrel{p}{\wedge} f_2 \cdots \stackrel{p}{\wedge}f_n = \big |A|_\infty \big( e_1^\star\stackrel{\infty}{\wedge} e_2^\star\cdots e_n^\star\stackrel{\infty}{\wedge} e_n^\star\big).
\end{equation}

Consequently, since   primal and  dual spaces are isomorphic, one can define for all $v_1,v_2, ...,v_n\in \Real^n$ the exterior product:
\begin{equation}
\big(v_1\stackrel{\infty}{\wedge} v_2\cdots \stackrel{\infty}{\wedge}v_n\big )=|v_1,v_2,...,v_n|_\infty\big(e_1\stackrel{\infty}{\wedge} e_2\cdots \stackrel{\infty}{\wedge}e_n\big ).
\end{equation}
Notice however, that though this definition extends as a limit case the usual definition of exterior product, it does not satisfy the the additivity property in each arguments with respect to the operation $\boxplus$.

\section{Semi-continuous Regularization and Limit of Hyperplanes
}\label{Reg}

\subsection{Semi-continuous Regularizations}
In the following, we say that a map $f:\Real^n\longrightarrow
\Real$ is a {\bf $\mathbb B$-form } if there exists some
$a\in\Real^n$ such that:
\begin{equation}f(x)=\bigboxplus_{i\in [n]}a_ix_i=\langle a,x\rangle_\infty.
\end{equation}
The function above is
depicted in Figure \ref{Rec}.2.

\medskip
\begin{center}{\scriptsize
\unitlength 0.5mm 
\linethickness{0.4pt}
\ifx\plotpoint\undefined\newsavebox{\plotpoint}\fi 
\begin{picture}(204.25,127.75)(-20,0)
\put(68.25,6.75){\vector(0,1){116.5}}
\put(0,55){\vector(1,0){141.75}}
\multiput(2.75,15.25)(.055722891566,.033734939759){2490}{\line(1,0){.055722891566}}
\put(31.25,77.25){\line(0,1){0}}
\put(184.25,77.25){\line(0,1){0}}
\put(51.25,83.25){\line(0,1){0}}
\put(204.25,83.25){\line(0,1){0}}
\put(112.5,81.5){\line(-1,0){83.75}}
\put(28.75,81.5){\line(0,-1){50.5}}
\put(96.75,72.25){\line(-1,0){52.75}}
\put(44,72.25){\line(0,-1){32}}
\put(128.75,91.5){\line(-1,0){115.75}}
\multiput(13,91.5)(.03125,-8.78125){8}{\line(0,-1){8.78125}}
\put(103.25,75.75){\line(0,-1){40.25}}
\put(103.25,35.5){\line(-1,0){67.5}}
\put(157.75,16.25){\line(0,1){0}}
\put(145.75,55.25){\makebox(0,0)[cc]{$x_1$}}
\put(69.25,127.75){\makebox(0,0)[cc]{$x_2$}}
\put(72.25,49.25){\makebox(0,0)[cc]{$0$}}
\put(21,97.5){\makebox(0,0)[cc]{$\langle a,x\rangle_\infty=c>0$}}
\put(88.25,30.75){\makebox(0,0)[cc]{$\langle a,x\rangle_\infty=d<0$}}
\put(147.25,109){\makebox(0,0)[cc]{$\langle a,x\rangle_\infty=0$}}
\put(68.25,0){\makebox(0,0)[cc]
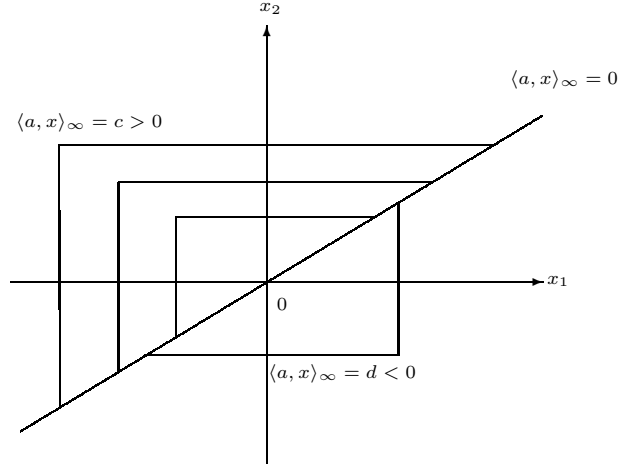
{{\bf Figure \ref{Reg}.1:} The level lines of the form $\langle a, \cdot\rangle_\infty$}}
\end{picture}}\end{center}
\medskip
These functions were used in \cite{b15,b17} to establish a separation theorem for $\mathbb B$-convex sets \cite{bh}.\footnote{A relaxed definition of $\mathbb B$-convexity was proposed in \cite{b15}: a subset $C$ of $\Real^n$ is $\mathbb B^\sharp$-convex if for all $x,y\in C$ and all $t\in [0,1]$, $x\boxplus ty\in C$.} All the points such that $\langle
a,x\rangle_{\infty}= 0$ are represented by the diagonal   line. In the
following, for all subsets $E$ of $\Real^n$ $\cl  (E)$ and $\inter
(E)$ respectively stand for the closure and the interior of $E$.


For all maps $f:\Real^n\longrightarrow \Real$ and
all real numbers $c$, the notation $[\,f\leq c\,]$ stands for the
set $f^{-1}(\,]-\infty, c]\,)$.  Similarly, $[\,f< c\,]$ stands
for  $f^{-1}(\,]-\infty, c[\,)$ and $[\,f\geq c\,]=[\,-f\leq
-c\,]$.

  For all $u,v\in \Real$, let us define the binary operation
\begin{equation*}\label{Bform}u\stackrel{-}{\smile} v=
\left\{\begin{matrix}
u &\hbox{ if } &|u|& > &|v|\\
\min \{u,v\}&\hbox{ if }&|u|&=&|v|\\
v& \hbox{ if }& |u|&<&|v|.\end{matrix}\right.\end{equation*} An
elementary calculus shows that $u\boxplus v=\frac{1}{2}\Big[u
\stackrel{-}{\smile} v-\big[(-u) \stackrel{-}{\smile}(-v)\big]
\Big]$.

 Similarly one can introduce a symmetrical binary operation defined for all $u,v\in \Real$ defined as:
\begin{equation*}\label{Bform}u\stackrel{+}{\smile} v=
\left\{\begin{matrix}
u &\hbox{ if } &|u|& > &|v|\\
\max \{u,v\}&\hbox{ if }&|u|&=&|v|\\
v& \hbox{ if }& |u|&<&|v|.\end{matrix}\right.\end{equation*}
Equivalently, one has: $u\stackrel{+}{\smile}
v=-\big[(-u)\stackrel{-}{\smile} (-v)\big]$. This means that $u\boxplus
v=\frac{1}{2}\Big[(u\stackrel{-}{\smile} v)+
(u\stackrel{+}{\smile} v)\Big]$.
Notice that the operations $\stackrel{-}{\smile}$ and $\stackrel{+}{\smile}$ are associative.
 Given $m$ elements $u_1, ..., u_m$ of $\Real$, not all of which are $0$, let $I_+$, respectively $I_-$,
 be the set of indices for which $0 < u_i$, respectively $u_i < 0$. We can then write
 $u_1\stackrel{-}{\smile}\cdots\stackrel{-}{\smile} u_m = (\stackrel{-}{\smile}_{i\in I_+}u_i) \stackrel{-}{\smile} (\stackrel{-}{\smile}_{i\in I_-}u_i)
 = (\max_{i\in I_+}u_i)\stackrel{-}{\smile}(\min_{i\in I_-}u_i)$ from which we have:
 \begin{equation}\label{manysmiles}
 u_1\stackrel{-}{\smile}\cdots\stackrel{-}{\smile} u_m =
\left\{
\begin{array}{lcc}
\max_{i\in I_+}u_i  &\hbox{if}   & I_- = \emptyset \hbox{ or } \max_{i\in I_-}\vert u_i\vert <  \max_{i\in I_+}u_i    \\
\min_{i\in I_-}u_i  &  \hbox{if} &  I_- = \emptyset \hbox{ or } \max_{i\in I_+}u_i <  \max_{i\in I_-}\vert u_i\vert \\
 \min_{i\in I_-}u_i  & \hbox{if}  &    \max_{i\in I_-}\vert u_i\vert  = \max_{i\in
 I_+}u_i.
\end{array}
\right.
 \end{equation}
We define a {\bf lower  ${\mathbb B}$-form }on ${\mathbb R}^n_+$
as a map $g: {\mathbb R}^n\to\Real$ such that for all $(x_1, ...,
x_n)\in{\mathbb R}^n_+$,
 \begin{equation}\label{eqcarBform}
 g(x_1, ..., x_n) = {\langle a,x\rangle}_\infty^{-} =a_1x_1 \stackrel{-}{\smile}\cdots \stackrel{-}{\smile} a_nx_n.
 \end{equation}
It was established in \cite{b15} that for all $c \in \Real$,
$g^{-1}\left(\,]-\infty,c]\right)=\left\{x\in \Real^n: g(x)\leq
c\right\}$ is closed. It follows that  a ${\mathbb B}$-form is
lower semi-continuous. It was established in \cite{bh3} that
$g^{-1}\left(\,]-\infty,c]\right)\cap \Real_+^n$ is a $\mathbb
B$-halfspace, that is a $\mathbb B$-convex subset of $\Real_+^n$
whose the complement in $\Real_+^n$ is also $\mathbb B$-convex.

Similarly, one can define an {\bf upper $\mathbb B$-form} as  a
map $h: {\mathbb R}^n\to\Real$  such that, for all $(x_1, ...,
x_n)\in{\mathbb R}^n$,
 \begin{equation}\label{eqcarBform}
 h(x_1, ..., x_n) = {\langle a,x\rangle}_\infty^{+} =a_1x_1 \stackrel{+}{\smile}\cdots \stackrel{+}{\smile} a_n x_n.
 \end{equation}
 For all $x\in \Real^n$, we clearly, have the following identities

 \begin{equation}
 \langle a,x\rangle_\infty^{+}=-\langle a,-x\rangle_\infty^{-}
 \text{ and }
 \langle a,x\rangle_\infty^{-}=-\langle a,-x\rangle_\infty^{+}.
 \end{equation}

The largest (smallest) lower (upper) semi-continuous minorant
(majorant) of a map $f$ is said to be the lower (upper)
semi-continuous regularization of $f$. In the next statements it
is shown that the lower (upper) $\mathbb B$-forms are the lower
(upper) semi-continuous regularized of the $\mathbb B$-forms.

\begin{prop} \cite{b15}\label{Reg}
\noindent Let $g$ be a lower $\mathbb B$-form defined by $ g(x_1,
..., x_n) = a_1x_1\stackrel{-}{\smile}\cdots\stackrel{-}{\smile}
a_nx_n, $ for some $a\in \Real^n$. Then $g$ is the lower
semi-continuous regularization of the map $x\mapsto \langle
a,x\rangle_{\infty}=\bigboxplus_{i\in [n]}a_ix_i$.
\end{prop}
The following corollary is then immediate.
\begin{cor} \cite{b15}\label{Reg2}
\noindent Let $h$ be an upper $\mathbb B$-form defined by $ h(x_1,
..., x_n) = a_1x_1\stackrel{+}{\smile}\cdots\stackrel{+}{\smile}
a_nx_n, $ for some $a\in \Real^n$. Then $h$ is the upper
semi-continuous regularization of the map $x\mapsto \langle
a,x\rangle_{\infty}=\bigboxplus_{i\in [n]}a_ix_i$.
\end{cor}

For the sake of simplicity let us denote for all $x\in \Real^n$:

\begin{equation} \bigsmileminus_{i\in I}x_i=x_1\stackrel{-}{\smile}\cdots\stackrel{-}{\smile} x_n \quad \text{and}\quad\bigsmileplus_{i\in I}x_i=x_1\stackrel{+}{\smile}\cdots\stackrel{+}{\smile} x_n.\end{equation}

 Let  $f:\Real^n\longrightarrow
\Real$ be a $\mathbb B$-form. Let $f^{-}$ and $f^{+}$ be
respectively the lower and upper semi-continuous regularized of
$f$ over $\Real^n$. It it shown in \cite{b17} that for all $c \in \Real$,
\begin{equation}\label{regmoinsplus}\cl{\big [f\leq c\,\big ]}={\big[ f^-
\leq  c \,\big]}\quad \text{ and }\quad \cl{\big [f\geq c\,\big ]}={\big[ f^+
\geq  c \,\big]}.\end{equation}
The following lemma is useful.

\begin{lem}\label{equal}For all dual $\mathbb B$-forms $f$ we have
$$\big[ f^-
\leq  0 \,\big]\cap \big[ f^+
\geq  0 \,\big]=\big[f^-+f^+=0\big].$$
\end{lem}
{\bf Proof:} Suppose that $x\in \big[ f^-
\leq  0 \,\big]\cap \big[ f^+
\geq  0 \,\big]$. If $f^-(x)=f^+(x)=0$, the inclusion is trivial. Suppose now that  $f^-(x)<0$ and  $f^+(x)>0$. There exists $a\in \Real^n$ such that
$f^-(x)=a_1x_1\stackrel{-}{\smile}\cdots \stackrel{-}{\smile}a_nx_n$ and $f^+(x)=a_1x_1\stackrel{+}{\smile}\cdots \stackrel{+}{\smile}a_nx_n$.
Hence there is some $i_-\in [n]$ such that $f^-(x)=a_{i_-}x_{i_-}<0$ and some  $i_+\in [n]$ such that $f^+(x)=a_{i_+}x_{i_+}>0$. However, since
by hypothesis this implies that $|a_{i_-}x_{i_-} |=|a_{i_+}x_{i_+}|$, we deduce that $ a_{i_-}x_{i_-}  =-a_{i_+}x_{i_+}$. Thus
$f^-(x)+f^+(x)=0$. Hence $\big[ f^-
\leq  0 \,\big]\cap\big[ f^+
\geq  0 \,\big]\subset \big[f^-+f^+=0\big]$. Conversely if $x\in \big[f^-+f^+=0\big]$, we have $f^-(x)f^+(x)\leq 0$, which implies the converse inclusion and ends the proof. $\Box$\\

In the remainder, it will be useful to consider the lower and upper semi-continuous determinant defined as:
\begin{equation}|A|_\infty^-=\bigsmileminus_{\sigma\in S_n} \big(\sgn(\sigma).
 \prod_{i\in [n]}a_{i,\sigma(i)}\big)\text{ and }|A|_\infty^+=\bigsmileplus_{\sigma\in S_n} \big(\sgn(\sigma).
 \prod_{i\in [n]}a_{i,\sigma(i)}\big).\end{equation}

\begin{center}{\scriptsize
\unitlength 0.29mm 
\linethickness{0.2pt}
\ifx\plotpoint\undefined\newsavebox{\plotpoint}\fi 
\begin{picture}(232.5,172.25)(0,0)
\put(107.75,6){\vector(0,1){160.5}}
\put(0,78.5){\vector(1,0){226}}
\multiput(3,19.5)(.05857988166,.03372781065){845}{\line(1,0){.05857988166}}
\put(52.5,48){\line(0,1){61.25}}
\put(52.5,109.25){\line(1,0){114.5}}
\multiput(167,109.25)(.0784552846,.0337398374){615}{\line(1,0){.0784552846}}
\put(50,132.75){\makebox(0,0)[cc]{$[f^+\geq c]$}}
\put(149.75,90.75){\makebox(0,0)[cc]{$[f^-\leq c]$}}
\put(225.25,138.5){\makebox(0,0)[cc]{$[f^-\leq c]\cap[f^+\geq c]$}}
\put(232.5,78){\makebox(0,0)[cc]{$x_1$}}
\put(106.5,172.25){\makebox(0,0)[cc]{$x_2$}}
\put(115.5,73.5){\makebox(0,0)[cc]{$0$}}
\put(81.25,-5){\makebox(0,0)[cc]
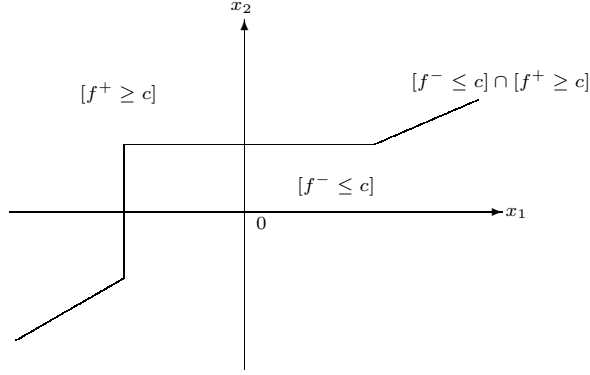
{{\bf Figure \ref{Reg}.1} Lower and Upper halfspaces.}}
\end{picture}}
\end{center}

\subsection{Kuratowski-Painlev\'e Limit of  Hyperplanes}\label{Limit}

This section is devoted to analyze the Kuratowski-Painlev\'e limit
of a sequence of half-spaces defined on the scalar field $(\Real,\stackrel{p}{+},\cdot)$.
The next result was established in \cite{b17}. These half-spaces are called $\varphi_p$-halfspaces.

\begin{prop} \label{limkur}Let $f$ be a
$\mathbb B$-form defined by $f(x)=\langle a,x \rangle_{\infty}$
for some $a\in \mathbb R^n\backslash\{0\}$. For any natural
number  $p$ let $f_{p}:\Real^n\longrightarrow \Real$ be a map
defined by $f_{p}(x)=\langle a_{p},x \rangle_{p}$ where
$\{a^{(p)}\}_{p\in \mathbb{N}}$ is a sequence  of
$\Real^n\backslash \{0\}$. If there exists a sequence
$\{c_{p}\}_{p\in \mathbb{N}}\subset \Real$  such that
$\lim_{q\longrightarrow \infty}(a^{(p) },c_{p })=(a,c)$, then:
$$Lim_{p\longrightarrow \infty}\big [f_p\leq c_p\big ]
=\cl{\big [f\leq c\big ]}=\big [f^-\leq c\big ]$$ and
$$Lim_{p\longrightarrow \infty}\big [f_p\geq c_p\big ]
=\cl{\big [f\geq c\big ]}=\big [f^+\geq c\big ].$$

\end{prop}

\begin{center}{\scriptsize
\unitlength 0.29mm 
\linethickness{0.3pt}
\ifx\plotpoint\undefined\newsavebox{\plotpoint}\fi 
\begin{picture}(232.5,172.25)(0,0)
\put(107.75,6){\vector(0,1){160.5}}
\put(0,78.5){\vector(1,0){226}}
\multiput(3,19.5)(.05857988166,.03372781065){845}{\line(1,0){.05857988166}}
\put(52.5,48){\line(0,1){61.25}}
\put(52.5,109.25){\line(1,0){114.5}}
\multiput(167,109.25)(.0784552846,.0337398374){615}{\line(1,0){.0784552846}}

\put(-15.75, 25.75){\makebox(0,0)[cc]{$[f_p\leq c]$}}

\put(149.75,90.75){\makebox(0,0)[cc]{$[f^-\leq c]$}}

\put(232.5,78){\makebox(0,0)[cc]{$x_1$}}
\put(106.5,172.25){\makebox(0,0)[cc]{$x_2$}}
\put(115.5,73.5){\makebox(0,0)[cc]{$0$}}
\put(81.25,0){\makebox(0,0)[cc]
{{\bf Figure \ref{Limit}} Limit of a sequence of $\varphi_p$-halfspaces.}}
\qbezier(3,23)(52.375,52.875)(52.25,78.25)
\qbezier(209.75,131.75)(166.75,109.5)(107.75,109.25)
\qbezier(107.75,109.25)(52.375,109)(52.5,78.75)
\qbezier(3,21)(52.25,43.375)(52.5,78.25)
\qbezier(52.5,78.25)(53.625,99)(58.25,102.75)
\qbezier(58.5,102.75)(58.875,108.875)(107.75,109.5)
\qbezier(107.75,109.5)(181.125,109.375)(210,129.75)
\qbezier(2.75,26.75)(48.625,58.875)(52,78.5)
\qbezier(52,78.5)(57.625,90.25)(62.75,95)
\qbezier(63,95)(66.375,101)(107.25,109)
\qbezier(107.25,109)(160,112.125)(207.75,133.75)
\end{picture}}
\end{center}

\bigskip
\bigskip

 In the following, one can go a bit further by showing that a sequence
 of
 $\phi_p$-hyperplanes defined for all $p\in \mathbb N$ as $[\langle
 a^{(p)},\cdot\rangle_p=c_p]$ has a Painlev\'e-Kuratowski limit.

\begin{prop}\label{limhyper} Let $f$ be a
$\mathbb B$-form defined by $f(x)=\langle a,x \rangle_{\infty}$
for some $a\in \mathbb R^n\backslash\{0\}$. For any natural
number  $p$ let $f_{p}:\Real^n\longrightarrow \Real$ be a map
defined by $f_{p}(x)=\langle a^{(p)},x \rangle_{p}$ where
$\{a^{(p)}\}_{p\in \mathbb{N}}$ is a sequence  of
$\Real^n\backslash \{0\}$. If there exists a sequence
$\{c_{p}\}_{p\in \mathbb{N}}\subset \Real$  such that
$\lim_{q\longrightarrow \infty}(a^{(p)},c_{p })=(a,c)$, then:
$$Lim_{p\longrightarrow \infty}\big [f_p=c_p\big ]
= \big [f^-\leq c\big ]\cap \big [f^+\geq c\big ].$$

\end{prop}
{\bf Proof:} By definition, for all $p$, we have $[f_p=c_p\big
]=[f_p\leq c_p\big ]\cap [f_p\geq c_p\big ]$. Hence, we have the
inclusion: \begin{align*}Ls_{p\longrightarrow \infty}[f_p=c_p\big
]&=Ls_{p\longrightarrow \infty}\Big([f_p\leq c_p\big ]\cap
[f_p\geq c_p\big ]\Big)\\ &\subset \Big(Ls_{p\longrightarrow
\infty} [f_p\leq c_p\big ] \Big)\cap \Big(Ls_{p\longrightarrow
\infty} [f_p\leq c_p\big ]\Big)\\&=\big [f^-\leq c\big ]\cap \big
[f^+\geq c\big ].\end{align*}
 In the following, we show that $\big [f^-\leq c\big ]\cap \big
[f^+\geq c\big ]\subset Li_{p\longrightarrow
\infty}[f_p=c_p\big].$ From Proposition \ref{limkur}, we have
$Li_{p\longrightarrow \infty}[f_p\leq c_p\big ]= \big [f^-\leq
c\big ]$ and $Li_{p\longrightarrow \infty}[f_p\geq c_p\big ]= \big
[f^+\geq c\big ]$. Suppose that $x\in \big [f^-\leq c\big ]\cap
\big [f^+\geq c\big ]$. This implies that there exist two
sequences $\{y^{(p)}\}_{p\in \mathbb N}$ and $\{z^{(p)}\}_{p\in
\mathbb N}$ respectively such that for any $p$, $y^{(p)}\in [f_p\leq c_p\big ]$
and  $z^{(p)}\in [f_p\geq c_p\big ]$ with
$x=\lim_{p\longrightarrow \infty}y^{(p)}=\lim_{p\longrightarrow
\infty}z^{(p)}$. For all $p$, the map $f_p$ is continuous.
Therefore, for all natural numbers $p$, there exists some
$\alpha_p\in [0,1]$ such that $f_p(\alpha_p
y^{(p)}+(1-\alpha_p)z^{(p)})=c_p$. Set $w^{(p)}=\alpha_p
y^{(p)}+(1-\alpha_p)^{(p)}$. We have for all natural numbers $p$
\begin{align*}\|x-w^{(p)}\|& =\|\alpha_p
(x-y^{(p)})+(1-\alpha_p)(x-z^{(p)})\|\\& \leq \alpha_p \|
x-y^{(p)}\| +(1-\alpha_p)\| x-z^{(p)}\|\\&\leq  \| x-y^{(p)}\| +\|
x-z^{(p)}\|.\end{align*} By hypothesis $\lim_{p\longrightarrow
\infty}\| x-y^{(p)}\|=\lim_{p\longrightarrow \infty}\|
x-z^{(p)}\|=0$. Thus $\lim_{p\longrightarrow \infty}\|
x-w^{(p)}\|=0$. Since $w_p\in [f_p= c_p\big ]$ for all $p$, we
deduce that $x\in Li_{p\longrightarrow \infty}[f_p=c_p\big]$. Consequently, $\big [f^-\leq c\big ]\cap \big [f^+\geq c\big
]\subset Li_{p\longrightarrow \infty}[f_p=c_p\big].$ Since we have
the sequence of inclusions
{ \begin{align*}Ls_{p\longrightarrow \infty}[f_p=c_p\big
]  \subset\big [f^-\leq c\big ]\cap \big [f^+\geq c\big
]\subset Li_{p\longrightarrow \infty}[f_p=c_p\big],\end{align*}} we deduce that
$$Lim_{p\longrightarrow \infty}[f_p=c_p\big
]=\big [f^-\leq c\big ]\cap \big [f^+\geq c\big ].\quad \Box$$

\subsection{ Limit Hyperplane Passing Through $n$ Points}

In this subsection we give the equation on a limit hyperplane passing though $n$ points.  Given $n$ points    $v_{1},..., v_n$ be $n$ in $\Real^n$, let $V$ be the $n\times n$ matrix whose each column is a vector $v_i$. If $|v_1,..., v_n|_p\not=0$ then let $H_p(V)$ denotes the $\varphi_p$-hyperplane passing trough $v_1,..., v_n$.

\begin{prop}Let  $v_1,..., v_n$ be $n$ points in $\Real^n$ and let $V$ be the $n\times n$ matrix whose each column $i$ is a vector $v_i$. Let $V_{(i)}$ be the matrix obtained from $V$ by replacing line $i$ with the transpose of the unit vector $1\!\!1_n$. Suppose that $|V|_\infty\not=0$. Then
$$Lim_{p\longrightarrow \infty}H_p(V)=\Big\{x\in \Real^n: \bigsmileminus_{i\in [n]} |V_{(i)}|_\infty x_i\leq |V|_\infty\leq \bigsmileplus_{i\in [n]} |V_{(i)}|_\infty x_i\Big\}.$$
\end{prop}
{\bf Proof:} First note that since $|V|_\infty\not=0$, there exists $p_0\in \mathbb N$ such that for all $p\geq p_0$, $|V|_p\not=0$. Therefore for all $p\geq p_0$, there exists a hyperplane $H_p(V)$ which contains $v_1,...,v_p$. Therefore, there exists some $a^{(p)}\in \Real^n$ and some  $c_p\in \Real$ such that
$$H_p(V)=\big [\langle a^{(p)}, \cdot\rangle_p=c\,\big ]. $$
Suppose that $x\in H_p(V)$. For all $i\in [n]$:
$$ \langle a^{(p)},v_i \stackrel{p}{-}x\rangle_p=c \stackrel{p}{-}c=0. $$
Let us denote $F_p(V)=\big [\langle a^{(p)}, \cdot\rangle_p=0\,\big ]$. Since $F_p(V)$ is a $n-1$-dimensional $\varphi_p$-subspace of $\Real^n$:
$$\big |v_1-x, v_2-x,\cdots, v_n-x\big|_p=0.$$

Let $V_{i,j}$ be the matrix obtained suppressing line $j$ and column $j$. It follows that:
\begin{align*}\big |v_1-x,& v_2-x,\cdots, v_n-x\big|_p\\&=\big |V\big |_p\stackrel{p}{-}\big |x, v_2,\cdots, v_n\big|_p\stackrel{p}{-}\big |v_1, x,v_3 , \cdots, v_n\big|_p\stackrel{p}{-}\big |v_1, v_2,\cdots,v_{n-1}, x\big|_p\\&=\big |V\big |_p \stackrel{p}{-} \stackrel{\varphi_p}{\sum_{j\in [n]}}\stackrel{\varphi_p}{\sum_{i\in [n]}}(-1)^{i+j}|V_{i,j}|_px_i =\big |V\big |_p \stackrel{p}{-} \stackrel{\varphi_p}{\sum_{i\in [n]}}\stackrel{\varphi_p}{\sum_{j\in [n]}}(-1)^{i+j}|V_{i,j}|_px_i \\&=\big |V\big |_p \stackrel{p}{-} \stackrel{\varphi_p}{\sum_{i\in [n]}}|V_{(i)}|_px_i=0 .\end{align*}
Therefore, we have:
$$H_p(V)=\Big\{ x\in \Real^n:\stackrel{\varphi_p}{\sum_{i\in [n]}}|V_{(i)}|_px_i=\big |V\big |_p\Big\}.$$
Since $\lim_{p\longrightarrow \infty}|V_{(i)}|_p=|V_{(i)}|_\infty$ and $\lim_{p\longrightarrow \infty}|V |_p=|V |_\infty$, we deduce the result from Proposition \ref{limhyper}. $\Box$\\

\medskip
A simple intuition is given in the case $n=2$ with two points. The hyperplane passing from two points $u$ and $v$ is a line. Let us denote
$D_p(u,v)$ the $\varphi_p$-line spanned by $u$ and $v$ in $\Real^2$. Every points $x=(x_1,x_2)\in \mathcal D_0(x,y)$ satisfy the relation:

\begin{equation}
| u {-} x^{}, v{-} x^{} | =\left |\begin{matrix}u_1- x_1^{}&v_1 -x_1\\u_2 -x_2^{}&v_2-x_2\end{matrix}\right|=0.
\end{equation}
Equivalently, we have
\begin{equation}
 (v_2-u_2)x_1^{}+(u_1-v_1)x_2^{}=|u,v|.
\end{equation}

For every points $z^{}\in   D_p(u,v)$ we have the relation:
\begin{equation*}
|u \stackrel{p}{-}x, v^{} \stackrel{p}{-}x|_p=\left |\begin{matrix}u_1^{}\stackrel{p}{-}x_1&v_1^{}\stackrel{p}{-}x_1\\u_2^{}\stackrel{p}{-}x_2&v_2^{}\stackrel{p}{-}x_2\end{matrix}\right|=0\iff (v_2\stackrel{p}{-}u_2)x_1^{ }\stackrel{p}{+}(u_1\stackrel{p}{-}v_1)x_2^{}=|u,v|_p.
\end{equation*}
We obtain that
\begin{align}&Lim_{p\longrightarrow \infty}D_p(u,v)\\&=\Big\{x\in \Real^2:\left| \begin{matrix} 1& 1\\u_2&v_2\end{matrix}\right|_\infty x_1\stackrel{-}{\smile}\left |\begin{matrix} u_1& v_1\\1&1\end{matrix} \right |_\infty x_2\leq |u,v|_\infty \leq \left| \begin{matrix} 1& 1\\u_2&v_2\end{matrix}\right|_\infty x_1\stackrel{+}{\smile}\left |\begin{matrix} u_1& v_1\\1&1\end{matrix} \right |_\infty x_2\Big\}.\notag \end{align}
Hence
\begin{align}&Lim_{p\longrightarrow \infty}D_p(u,v)\\&=
\Big\{x\in \Real^2: (  v_2\boxminus u_2)x_1\stackrel{-}{\smile}( u_1\boxminus v_1 ) x_2\leq |u,v|_\infty \leq (  v_2\boxminus u_2)x_1\stackrel{+}{\smile}( v_1\boxminus u_1 ) x_2\Big\}.\notag \end{align}

\begin{expl}Suppose that $n=3$ and that $v_1=(1,0,-3)$, $v_2=(2,-1,1)$, $v_3=(4,1,2)$. We have
$$V=\begin{pmatrix}1&2&4\\0&-1&1\\-3&1&2\end{pmatrix}.$$ Thus $$ V_{(1)}=\begin{pmatrix}1&1&1\\0&-1&1\\-3&1&2\\\end{pmatrix}, V_{(2)}=\begin{pmatrix}1&2&4\\1&1&1\\-3&1&2\\\end{pmatrix} \text{ and }\; V_{(3)}=\begin{pmatrix}1&2&4\\0&-1&1\\1&1&1\\\end{pmatrix}. $$
Hence $|V|_\infty=1.(-1). 2\boxplus 2.1.(-3)\boxplus 0.1.4 \boxplus (-4).(-1).(-3)\boxplus (-1).1.1 \boxplus 0.2.2=-12$; $|V_{(1)}|_\infty=1.(-1). 2\boxplus 1.1.(-3)\boxplus 0.1.1 \boxplus (-1).(-1).(-3)\boxplus (-1).1.1 \boxplus 0.1.2=-3$; $|V_{(2)}|_\infty=1.1. 2\boxplus 2.1.(-3)\boxplus 1.1.4 \boxplus (-4).1.(-3)\boxplus (-1).1.1 \boxplus (-1).2.2=12$; $|V_{(3)}|_\infty=1.(-1). 1\boxplus 2.1.1\boxplus 0.1.4 \boxplus (-4).(-1).1\boxplus (-1).1.1 \boxplus 0.2.1=4$.
$$H_\infty(V)=\big\{(x_1,x_2,x_3): (-3)x_1\stackrel{-}{\smile} 12 x_2 \stackrel{-}{\smile} 4 x_3\leq -12\leq (-3)x_1\stackrel{+}{\smile} 12 x_2 \stackrel{+}{\smile} 4 x_3\big \}.$$
It is easy to check that $v_1,v_2,v_3\in H_\infty(V)$.
\end{expl}

 \section{Limit Systems of Equations}\label{SystEq}

 \subsection{Limit  Solutions of a Sequence of Systems of Equations}

   For all $\varphi_p$-linear endomorphisms $f:x\mapsto A\stackrel{p}{\cdot}x$, where $A\in \mathcal M_{n}(\Real)$, we consider a sequence of $\varphi_p$-linear systems of the form $A\stackrel{p}{\cdot}x=b$.

   \begin{prop}\label{semisyst}Let $A\in \mathcal M_n(\Real)$ be a square matrix. If $|A|_\infty\not=0$ then there exists a uniqueness  $x^\star =\sum_{i\in [n]}\frac{|A^{(i)}|_\infty}{|A|_\infty}e_i$ such that

   $$\{x^\star\}=Lim_{p\longrightarrow
  \infty}\big\{x\in \Real^n:A\stackrel{p}{\cdot}x=b \big\}.$$
  Conversely, if there is some $x^\star\in \Real^n$ such that:
  $$ x^\star \in Ls_{p\longrightarrow
  \infty}\big\{x\in \Real^n:A\stackrel{p}{\cdot}x=b \big\}$$ then $|A|_\infty\not=0$ and $\{x^\star\}=Lim_{p\longrightarrow
  \infty}\big\{x\in \Real^n:A\stackrel{p}{\cdot}x=b \big\}.$

\end{prop}
{\bf Proof:} If $|A|_\infty\not=0$ then there is some $p_0$ such
that for all $p\geq p_0$, $|A|_p\not=0$. Thus for all $p\geq p_0$,
$x^{(p)}=\sum_{i\in [n]}\frac{|A^{(i)}|_p}{|A|_p}e_i$ is
solution of the system $A\stackrel{p}{\cdot}x=b$ and therefore
$x^{(p)}\in \big\{x\in \Real^n:A\stackrel{p}{\cdot}x=b \big\}$.
However, $x^\star =\sum_{i\in
[n]}\frac{|A^{(i)}|_\infty}{|A|_\infty}e_i=\lim_{p\longrightarrow
\infty}x^{(p)}$. Thus $x^\star\in Li_{p\longrightarrow
  \infty}\big\{x\in \Real^n:A\stackrel{p}{\cdot}x=b \big\}$. Moreover,  for all $p\geq p_0$, since  $|A|_p\not=0$ we have $\{x^{(p)}\}=\big\{x\in \Real^n:A\stackrel{p}{\cdot}x=b \big\}$. Consequently $x^\star$ is the uniqueness solution. This implies that
  $x^\star\in Ls_{p\longrightarrow
  \infty}\big\{x\in \Real^n:A\stackrel{p}{\cdot}x=b \big\}$. Moreover, for all increasing sequence of natural numbers $\{p_k\}_{k\in \mathbb N}$, $x^{(p_k)}=\sum_{i\in [n]}\frac{|A^{(i)}|_{p_k}}{|A|_{p_k}}e_i$ is the uniqueness solution of the system of the form $A\stackrel{p_k}{\cdot}x=b$. Hence $Li_{p\longrightarrow \infty}\big\{x\in \Real^n:A\stackrel{p}{\cdot}x=b \big\}=Ls_{p\longrightarrow \infty}\big\{x\in \Real^n:A\stackrel{p}{\cdot}x=b \big\}=\{x^\star\}$ which ends the first part of the statement.

   To complete the proof, suppose that $\{x^\star\}=Ls_{p\longrightarrow
  \infty}\big\{x\in \Real^n:A\stackrel{p}{\cdot}x=b \big\}$ with $|A|_\infty=0$ and let us show a contradiction. This implies that for any $p\in \mathbb N$  we have $|A|_p=0$. Thus, for any $p$, the system $\{x\in \Real^n:A\stackrel{p}{\cdot}x=b \big\}$ has either an infinity of solutions or is an empty set. If, for all $p\in \mathbb N$, it is an empty set then the upper limit of the sequence of solution sets is empty. Suppose that this is not the case and let us show a contradiction.  Suppose that $x^\star\in Ls_{p\longrightarrow
  \infty}\big\{x\in \Real^n:A\stackrel{p}{\cdot}x=b \big\}$. In such case there exists  a subsequence $\{p_k\}_{k\in\mathbb N}$ such that $x^\star=\lim_{k\longrightarrow \infty}x^{(p_k)}$ where for all $k$, $x^{(p_k)}\in \{ x\in \Real^n :A\stackrel{p_k}{\cdot}x=b\big\}$ that is a $\varphi_{p_k}$-affine subspace that contains an infinity of points. For any $k$ let us consider the ball $B_\infty(x^{(p_k)},1]$ of center $x^{(p_k)}$ and of radius $1$. Since $\{ x\in \Real^n :A\stackrel{p_k}{\cdot}x=b\big\}$ is a $\varphi_{p_k}$ affine subspace of $\Real^n$, for all $k$ there exits a vector $v^{(p_k)}\not=0$ such that $A\stackrel{p_k}{\cdot}v^{(p_k)}=0$. This implies that:
  $$\{x^{(p_k)}\stackrel{p_k}{+}\delta v{(p_k)}:\delta \in \Real\}\subset \{x\in \Real^n:
A\stackrel{p_k}{\cdot}x=b\big\} .$$ Let $\delta_{p_k}=\sup \{\delta: x^{(p_k)}\stackrel{p_k}{+}\delta v^{(p_k)}\in B_\infty(x^{(p_k)},1]\}$. Since the map $\delta\mapsto x^{(p_k)}\stackrel{p_k}{+}\delta v^{(p_k)}$ is a continuous vector valued function,  $y^{(p_k)}=x^{(p_k)}\stackrel{p_k}{+}\delta_{p_k} v^{(p_k)}\in C_\infty(x^{(p_k)},1)$ which implies that $d_\infty (x^{(p_k)},y^{(p_k)})=1$. Now since $\{x^{(p_k)}\}_{k\in \mathbb N}$ converges to $x^\star$, There exists some $d>0$ and $k_d\in \mathbb N$ such that for all $k\geq k_d$, $x^{(p_k)},y^{(p_k)}\in B_\infty(x^\star, d\,]$. Since   $B_\infty(x^\star, d\,]$ is a compact subset of $\Real^n$ one can extract a sequence $\{y^{(p_{k_r})}\}_{r\in \mathbb N}$ which converges to some $y^\star\in Ls_{p\longrightarrow
  \infty}\big\{x\in \Real^n:A\stackrel{p}{\cdot}x=b \big\}$. However, for all $r\in \mathbb N$, $d_\infty (x^{(p_{k_r})},y^{(p_{k_r})})=1$, and we deduce that $d_\infty (x^{ \star},y^{\star})=1$. This implies that $x^\star\not=y^\star$ which contradicts the unicity. Consequently if the upper limit of the sequence of solution sets has a uniqueness element, then $|A|_\infty\not=0$.  $\Box$\\

  In the following, for all matrices $A\in \mathcal M_{n,l}(\Real)$ and $B\in \mathcal M_{l,m}(\Real)$ let us define  the   product:
\begin{equation}
  A{\boxtimes}B=\Big(\bigboxplus_{k\in [l]}a_{i,k}b_{k,m}\Big)_{\substack{i\in [n]\\ j\in [m]}}.
  \end{equation}
  The lower and upper semi-continuous regularized products are respectively defined as:
  \begin{equation}
  A\stackrel{-}{\boxtimes}B=\Big(\bigsmileminus_{k\in [l]}a_{i,k}b_{k,m}\Big)_{\substack{i\in [n]\\ j\in [m]}}\quad \text{ and }\quad
  A\stackrel{+}{\boxtimes}B=\Big(\bigsmileplus_{k\in [l]}a_{i,k}b_{k,m}\Big)_{\substack{i\in [n]\\ j\in [m]}}.
  \end{equation}
By construction, it follows that for all vectors $x\in \Real^n$, the matrix-vector products derived from $\stackrel{-}{\boxtimes}$ and $\stackrel{+}{\boxtimes}$ are
  defined by:
  \begin{equation}
  A \stackrel{-}{\boxtimes} x= \begin{pmatrix}
{\langle a_1,x\rangle}_\infty^{-} \\

\vdots  \\
{\langle a_n,x\rangle}_\infty^{-} \end{pmatrix}
\text{  and }
  A \stackrel{+}{\boxtimes} x= \begin{pmatrix}
{\langle a_1,x\rangle}_\infty^{+} \\
\vdots  \\
{\langle a_n,x\rangle}_\infty^{+} \end{pmatrix}.
  \end{equation}
The next result is an immediate consequence.
\begin{prop}Let $A\in \mathcal M_n(\Real)$ be a square matrix. Suppose that
$x^\sharp\in Ls_{p\longrightarrow \infty}\{x\in \Real^n; A\stackrel{p}{\cdot}x=b\}$. Then $x^\sharp$ is solution of the system:
\begin{equation}\label{IS}\Bigg\{\begin{matrix}& A \stackrel{-}{\boxtimes} x&\leq& b&\\& A \stackrel{+}{\boxtimes} x&\geq& b,& x\in \Real^n.\end{matrix}\end{equation}
Moreover, if $|A|_\infty\not=0$ then  $x^\star =\sum_{i\in [n]}\frac{|A^{(i)}|_\infty}{|A|_\infty}e_i$ is solution of system \eqref{IS}.\end{prop}
{\bf Proof:} From Proposition \ref{limhyper}, for all $i\in [n]$:
  \begin{align*}Ls_{p\longrightarrow \infty}\{x: \langle a_i,x\rangle_p=b_i\}&=Lim_{p\longrightarrow \infty}\{x: \langle a_i,x\rangle_p=b_i\}\\&=[{x: \langle a_i,x\rangle}_\infty^{-}\leq b_i] \cap
[{x: \langle a_i,x\rangle}_\infty^{+}\geq b_i]. \end{align*}
From
\eqref{inclus}, we deduce that:
\begin{align*}Ls_{p\longrightarrow
  \infty}\big\{x\in \Real^n:A\stackrel{p}{\cdot}x=b \big\}&=Ls_{p\longrightarrow
  \infty}\bigcap_{i\in [n]}[\langle a_i,\cdot \rangle_p=b_i ]\\
  &\subset\bigcap_{i\in [n]}Ls_{p\longrightarrow
  \infty}[\langle a_i,\cdot \rangle_p=b_i ]\\
  &\subset \bigcap_{i\in [n]}\Big( [{ \langle a_i,\cdot \rangle}_\infty^{-}\leq b_i] \cap
[{  \langle a_i,\cdot\rangle}_\infty^{+}\geq b_i]\Big), \end{align*}
which implies that if $x^\sharp\in Ls_{p\longrightarrow
  \infty}\big\{x\in \Real^n:A\stackrel{p}{\cdot}x=b \big\}$ then it satisfies the system \eqref{IS}.
If $|A|_\infty\not=0$, from Proposition \ref{semisyst} $\{x^\star\} = Lim_{p\longrightarrow
  \infty}\big\{x\in \Real^n:A\stackrel{p}{\cdot}x=b \big\}$ and this implies that $x^\star$ satisfies system \eqref{IS}.
$\Box$\\

Since it contains any element of the upper limit set
$Ls_{p\longrightarrow \infty}\{x\in \Real^n:A\stackrel{p}{\cdot}x=b\}$, the system \eqref{IS} is called a {\bf limit system}.

 \begin{expl} Let us consider the matrix
 $$A=\begin{pmatrix}-1&1\\1&1\end{pmatrix}\quad$$
with $a_1=(-1,1)$, $a_2=(1,1)$ and suppose that $b_1=2$, $b_2=3$. Now, let us consider the matrices:
$$ A^{(1)}=\begin{pmatrix}2&1\\3&1\end{pmatrix}\quad
A^{(2)}=\begin{pmatrix}-1&2\\1&3\end{pmatrix}.$$
 We have $\left|A\right|_\infty
 =((-1)\cdot 1)\boxplus ((- 1)\cdot
 4)=-1$; $\left|A^{(1)}\right|_\infty=(2\cdot 1\boxplus (- 3\cdot
1))=-3$;
$\left|A^{(2)}\right|_\infty=((-1)\cdot
3\boxplus( (- 2)\cdot
 1)=-3$. We obtain the solutions: $x_1^\star=\frac{-3}{-1}=3\quad x_2^\star=\frac{-3}{-1}={3}.$ One can then check that:
$$\begin{pmatrix}-1&1\\1&1\end{pmatrix}\stackrel{-}{\boxtimes}\begin{pmatrix}3\\3\end{pmatrix}=\begin{pmatrix}(-3)\stackrel{-}{\smile} 3\\3\stackrel{-}{\smile} 3\end{pmatrix}=\begin{pmatrix}-3\\3\end{pmatrix}\leq \begin{pmatrix}2\\3\end{pmatrix}$$
and
$$\begin{pmatrix}-1&1\\1&1\end{pmatrix}\stackrel{+}{\boxtimes}\begin{pmatrix}3\\ 3\end{pmatrix}=\begin{pmatrix}(-3)\stackrel{+}{\smile} 3\\3\stackrel{+}{\smile} 3\end{pmatrix}=\begin{pmatrix}3\\3\end{pmatrix}\geq \begin{pmatrix}2\\3\end{pmatrix}.$$
Therefore $x^\star=(3,3)$ is a solution of the limit system. This example is depicted in Figure \ref{SystEq}.2.

\begin{center}{\scriptsize

\unitlength 0.4mm 
\linethickness{0.4pt}
\ifx\plotpoint\undefined\newsavebox{\plotpoint}\fi 
 \ifx\plotpoint\undefined\newsavebox{\plotpoint}\fi 
\begin{picture}(208.25,173.75)(0,0)
\put(103.25,18.25){\vector(0,1){147}}
\put(24.75,95.25){\vector(1,0){177.5}}
\multiput(22,22)(.03774089936,.03372591006){934}{\line(1,0){.03774089936}}
\put(57.25,53.5){\line(0,1){74.5}}
\put(57.25,128){\line(1,0){74.75}}
\multiput(132,128)(.03476123596,.03370786517){712}{\line(1,0){.03476123596}}
\put(185.75,134.25){\line(0,-1){.25}}
\put(146,141.25){\circle*{1.118}}
\put(89.75,-10){\makebox(0,0)[cc]
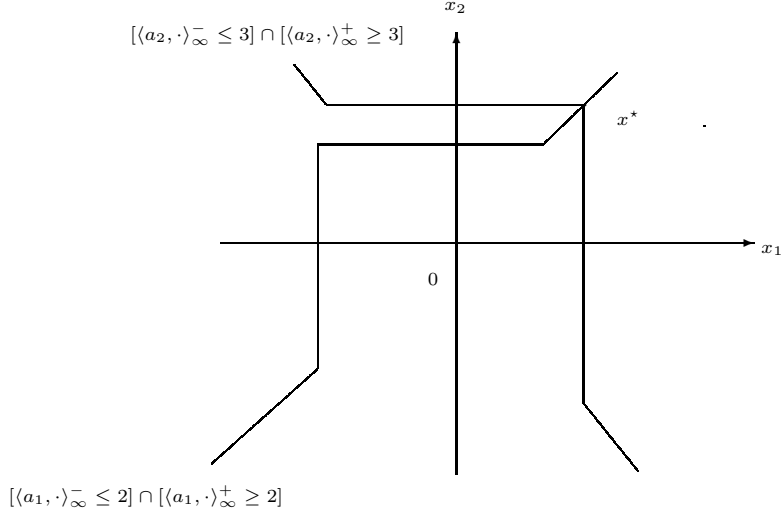
{{{\bf Figure \ref{SystEq}.2} Example of a two dimensional Limit System.}}}}
\put(208.25,93){\makebox(0,0)[cc]{$x_1$}}
\put(103,173.75){\makebox(0,0)[cc]{$x_2$}}
\put(95.5,83.25){\makebox(0,0)[cc]{$0$}}
\put(160.5,137){\makebox(0,0)[]{$x^\star$}}
\multiput(49.25,154.75)(.0336990596,-.0423197492){319}{\line(0,-1){.0423197492}}
\put(60,141.25){\line(1,0){85.5}}
\put(145.5,141.25){\line(0,-1){99.25}}
\multiput(145.5,42)(.0337338262,-.042051756){541}{\line(0,-1){.042051756}}
\put(40.5,164.75){\makebox(0,0)[cc]
{{$[\langle a_2,\cdot\rangle_\infty^-\leq 3]\cap[\langle a_2,\cdot\rangle_\infty^+\geq 3]$}}}
\put(0,11.25){\makebox(0,0)[cc]
{{$[\langle a_1,\cdot\rangle_\infty^-\leq 2]\cap[\langle a_1,\cdot\rangle_\infty^+\geq 2]$}}
\end{picture}
}
\end{center}

 \end{expl}

 \bigskip

\begin{expl}
Let us consider the matrix: $$A=\begin{pmatrix}3&-1&3\\2&-4&1\\-4&5&3\end{pmatrix}$$
with $a_1=
(3,1,-3)$, $a_2= (2,-4,1)$, $a_3=(-4,5,3)$, $b_1=6$, $b_2=8$, $b_3=4$.

The limit system is:

\begin{equation}\label{eq33}\left\{\begin{matrix}\begin{pmatrix}3&-1&3\\2&-4&1\\-4&5&3\end{pmatrix} \stackrel{-}{\boxtimes}\begin{pmatrix}x_1\\x_2\\x_3\end{pmatrix}\leq \begin{pmatrix}8\\4\\6\end{pmatrix}&\\\begin{pmatrix}3&-1&3\\2&-4&1\\-4&5&3\end{pmatrix} \stackrel{+}{\boxtimes}\begin{pmatrix}x_1\\x_2\\x_3\end{pmatrix}\geq \begin{pmatrix}8\\4\\6\end{pmatrix}& x\in \Real^3.\end{matrix}\right.\end{equation}

Now, let us consider the matrices:

 $$A^{(1)}=\begin{pmatrix}6&-1&3\\8&-4&1\\ 4&5&3\end{pmatrix};\quad A^{(2)}=\begin{pmatrix}3&6&3\\2&8&1\\-4&4&3\end{pmatrix}; A^{(3)}=\begin{pmatrix}3&-1&6\\2&-4&8\\-4&5&4\end{pmatrix}.$$
We have
$|A|_\infty=(-36)\boxplus4\boxplus 30\boxplus (-48) \boxplus (-15)\boxplus 6=-48$;
$|A^{(1)}1|_\infty=(-72)\boxplus (-4) \boxplus 120\boxplus 48 \boxplus 24\boxplus (-30)=120$;
$|A^{(2)}|_\infty= 72 \boxplus (-24) \boxplus 24\boxplus 96 \boxplus (-12)\boxplus (-36)=96$;
$|A^{(3)}|_\infty= (-48) \boxplus 32 \boxplus +60\boxplus (-96) \boxplus 8\boxplus (-120)=-120$.
We obtain that
$$x_1^\star =\frac{120}{-48}=-\frac{5}{2},\quad x_2^\star= \frac{96}{-48}={-2}\quad \text{ and }  \quad x_3^\star =\frac{5}{2}.$$
 Let us check that $x^\star=(-\frac{5}{2},-2,-\frac{5}{2})$ satisfies the system of equations \eqref{eq33}. We have:

 $$\left\{\begin{matrix}\begin{pmatrix}3&-1&3\\2&-4&1\\-4&5&3\end{pmatrix} \stackrel{-}{\boxtimes}\begin{pmatrix}-\frac{5}{2}\\-2\\\frac{5}{2}\end{pmatrix}=\begin{pmatrix}-\frac{15}{2}\\8\\-10\end{pmatrix}\leq \begin{pmatrix}6\\8\\4\end{pmatrix}&\\\begin{pmatrix}3&-1&3\\2&-4&1\\-4&5&3\end{pmatrix} \stackrel{+}{\boxtimes}\begin{pmatrix}-\frac{5}{2}\\-2\\\frac{5}{2}\end{pmatrix}=\begin{pmatrix} \frac{15}{2}\\8\\ 10\end{pmatrix}\geq \begin{pmatrix}6\\8\\4\end{pmatrix}& x\in \Real^3.\end{matrix}\right.$$
Thus $x^\star=(-\frac{5}{2},-2,-\frac{5}{2})$ satisfies the system \eqref{eq33}.\end{expl}
In the following, we say that a solution $x^\star$ of the limit
system is {\bf regular} if for all $i\in [m]$, $\langle
a_{i},x^\star \rangle_{\infty}=\langle a_{i},x^\star
\rangle_{\infty}^-=\langle a_{i},x^\star \rangle_{\infty}^+$. This
implies that $x^\star$ is also solution of the equation
\begin{equation}A\boxtimes x=b.\end{equation}
Equivalently, this means that:

\begin{equation}\left\{
\begin{matrix}
\bigboxplus_{j\in [n]}a_{1,j}x_j&=b_{1}\\
\vdots&\vdots\\
\bigboxplus_{j\in [n]}a_{m,j}x_j&=b_{m}.\\
\end{matrix}\right.
\end{equation}

 \begin{center}
 {\scriptsize \unitlength 0.4mm 
\linethickness{0.4pt}
\ifx\plotpoint\undefined\newsavebox{\plotpoint}\fi 
\begin{picture}(183.5,172.25)(0,0)
\put(78.5,16.75){\vector(0,1){147}}
\put(0,88.75){\vector(1,0){177.5}}
\put(161,132.75){\line(0,-1){.25}}
\put(183.5,89.5){\makebox(0,0)[cc]{$x_1$}}
\put(78.25,172.25){\makebox(0,0)[cc]{$x_2$}}
\put(70.75,81.75){\makebox(0,0)[cc]{$0$}}
\put(128.5,60.75){\makebox(0,0)[]{$x^\star$}}
\multiput(24.5,153.25)(.0336990596,-.0423197492){319}{\line(0,-1){.0423197492}}
\put(35.25,139.75){\line(1,0){85.5}}
\put(120.75,139.75){\line(0,-1){99.25}}
\multiput(120.75,40.5)(.0337338262,-.042051756){541}{\line(0,-1){.042051756}}
\put(24.25,53.75){\line(1,0){109.75}}
\multiput(134.25,53.75)(.0598591549,-.0337022133){497}{\line(1,0){.0598591549}}
\put(24,53.5){\line(0,1){68}}
\multiput(24,121.5)(-.0536144578,.0337349398){415}{\line(-1,0){.0536144578}}
\put(76.75,0){\makebox(0,0)[cc]
{{\bf
Figure \ref{SystEq}.3} Regular Solutions of a Limit System.}}
\put(120.5,54){\circle*{1.5}}
\end{picture}} \end{center}

\subsection{Positive Solutions of Positive Systems of Maximum Equations }

In the following, we consider a theorem established by Kaykobad
\cite{kayko} that gives a necessary condition for the existence of
a positive solution to a positive invertible linear system.

\begin{thm} \label{condpos} Suppose that
$A=(a_{i,j})_{i,j\in [n]}\in \mathcal M_n(\Real) )$ is a square
matrix such that for all $i,j$ $a_{i,j}\geq 0$  and $ a_{i,i} > 0$
for all $i\in [n]$. Suppose moreover that $b\in \Real_{++}^n$. If
for all $i\in [n]$
$$b_i > \sum_{j\in [n]\backslash \{i\}}a_{i,j}\frac{b_j}{a_{j,j}} $$ then $A$ is invertible
 and $A^{-1}b  \in \Real_{++}^n $.
\end{thm}

In the following this result is extended to a
$\varphi_p$-endomorphism.

\begin{lem}  Suppose that
$A=(a_{i,j})_{i,j\in [n]}\in \mathcal M_n(\Real)  $ is a square
matrix such that for all $i,j$ $a_{i,j}\geq 0$. Suppose that there exists a permutation $\sigma:[n]\longrightarrow [n]$ such that $ a_{i,\sigma(i)} > 0$
for all $i\in [n]$. Suppose moreover that $b\in \Real_{++}^n$. If
for all $i\in [n]$
$$b_i > \Big(\sum_{j\in [n]\backslash \{i\}}{(a_{i,\sigma(j)})}^{2p+1}\frac{{(b_j)}^{2p+1}}{{(a_{j,\sigma(j)})}^{2p+1}}\Big)^{\frac{1}{2p+1}} $$ then $A$ is $\varphi_p$-invertible and there is a solution $x^{(p)}\in \Real_{++}^n$ to the equation $A\stackrel{p}{\cdot}x=b$.
\end{lem}
{\bf Proof:} Let $\bar A=(\bar a_{i,j})_{i,j\in [n]}$ be the $n\times n$ matrix defined by $ \bar a_{i,j}= a_{i,\sigma(j)}$. The system $\bar A\stackrel{p}{\cdot}x=  b$ is equivalent to
$\Phi_p(\bar A)u=\phi_p(  b)$ setting $u=\phi_p(x)$. Since $ a_{i,\sigma(i)} > 0$
for all $i$, we deduce that for all $i$, $\bar a_{i,i}>0$. Since by definition
$\Phi_p(\bar A)=({{\bar a_{i,j}}}^{2p+1})_{i,j\in [n]}$ and
$\phi_p(b)=({  b_1}^{2p+1},...,  {b_n}^{2p+1})$, it follows from Theorem
\ref{condpos} that this system has a positive solution if:
$$( b_i)^{ {2p+1}} > \sum_{j\in [n]\backslash \{i\}}{(\bar a_{i,j})}^{2p+1}\frac{{(  b_j)}^{2p+1}}{{(\bar a_{j,j})}^{2p+1}}.  $$
Equivalently, we deduce that the system   $ A\stackrel{p}{\cdot}x=  b$ has a solution if $$( b_i)^{ {2p+1}} > \sum_{j\in [n]\backslash \{i\}}{(a_{i,\sigma(j)})}^{2p+1}\frac{{(b_j)}^{2p+1}}{{(a_{j,\sigma(j)})}^{2p+1}}  $$ which ends the proof. $\Box$\\

First, we consider systems of max-equations, that is, systems of
the form
\begin{equation}\label{MaxEq}\left\{
\begin{matrix}
\max\{a_{1, 1}x_{1},\ldots, a_{1, n}x_{n}\}&=b_{1}\\
\vdots&\vdots\\
\max\{a_{m, 1}x_{1},\ldots, a_{m, n}x_{n}\}&=b_{m}\\
\end{matrix}\right.\end{equation}
where $a_i = (a_{i, 1},\ldots, a_{i, n})\in\Real_{+}^{n}, i = 1,
\ldots, m$, $b = (b_1, \ldots, b_m)\in\Real_{+}^{m}$ and the
solution $(x_1, \ldots, x_n)$ is to be found in $\Real_{+}^{n}$.
Notice that if $b_i = 0$ then we have to take $x_j = 0$ for each
$j$ such that $a_{i, j} > 0$, and, as far as equation $i$ is
concerned, the other values  $x_l$ are irrelevant; equation $i$
can therefore be removed from the system and the number of
variables decreases. In other words, we can assume that $b_i > 0$
for all $i$. In the remainder these types of systems will called system of maximum-equations.
We can assume that for all $j$ there is at least one index $i$
such that $a_{i, j} > 0$;  let $\eta(j) = \{i : a_{i, j} > 0
\}$ and
\begin{equation}x^\star = \sum_{i\in [n]}\big(\min_{i\in \eta (j)} \frac{b_i}{a_{i, j}} \big)e_i.\end{equation}
From \cite{bh}, the  system of maximum equations  \eqref{MaxEq} has some solution, then $x^\star$ is a solution and, for any solution $x$ one has $x\leq x^\star.$ This condition is equivalent to the following.

\begin{lem}\label{eqforms}Let $A\in \mathcal M_n(\Real)$ be a square matrix such that $a_{i,j}\geq 0$ for all $i,j\in [n]$.
For all $i,j\in [n]$, let us denote $\mu(i) = \{j : a_{i, j} > 0
\}$ and $\eta(j) = \{i : a_{i, j} > 0
\}$ and  assume that   $ \eta(j)$ and $ \mu(i)$ are nonempty. Suppose moreover that $b\in \Real_{++}^n$. The system of maximum equations  \eqref{MaxEq} has a solution in $\Real_+^n$ if and only if there exists a permutation $\sigma:
[n]\longrightarrow [n]$ such that for all $i\in [n]$
$$b_i \geq \max_{j\in [n]\backslash \{i\}}{ a_{i,\sigma(j)} } \frac{{ b_j }}{{ a_{j,\sigma(j)} } }. $$
Moreover, this solution is uniqueness if and only if for all $i\in [n]$
$$b_i > \max_{j\in [n]\backslash \{i\}}{ a_{i,\sigma(j)} } \frac{{ b_j }}{{ a_{j,\sigma(j)} } }. $$
\end{lem}
{\bf Proof:}  The system \eqref{MaxEq} has a solution if and only if the point $x^\star=\sum_{j\in [n]} \big(\min_{i\in \eta (j)} \frac{b_i}{a_{i, j}}\big) e_j$ is a solution. Suppose  that $x^\star$ is solution. Let us assume that there exists $j\in [n]$ such that for all
$k$ and all $i\not=j$, we have  $\frac{a_{j,k}}{b_j}< \frac{a_{i,k}}{b_{i}}$
and let us show a contradiction. This implies that for all $k$, $ \frac{a_{j,k}}{b_j} <  \max_{i\in [n]}\big(\frac{a_{i,k}}{b_{i}} \big) $. Therefore for all $k\in \mu(j)$, $ \min_{i\in \eta(k)}\big(\frac{b_{i}}{a_{i,k}} \big)<  \frac{b_j}{a_{j,k}}$. Set $x^{(j)}=\sum_{k\in \mu(j)}\frac{b_j}{a_{j,k}}e_k$. Hence $\max_{k\in [n]}(a_{j,k}x_k^\star)<
\max_{k\in \mu(j)}(a_{j,k}x_k^{(j)})=b_j$. However, since $x^\star$ is solution of system \eqref{MaxEq}, this is a contradiction. Hence, for all $j$, there exists $\sigma(j)$ such that for all $i\not=j$, we have
$$\frac{a_{j,\sigma(j)}}{b_j}\geq \frac{a_{i,\sigma(j)}}{b_{i}}.$$ Since for all $j$ we have $\eta(j)\not=\emptyset$, we deduce  that, for all $j$,  $a_{j,\sigma(j)}>0$. Therefore, this is equivalent to the condition $b_i \geq { a_{i,\sigma(j)} } \frac{{ b_j }}{{ a_{j,\sigma(j)} }}$ for all $j\not=i$. Consequently, we deduce that
$$ b_i \geq \max_{j\in [n]\backslash \{i\}}{ a_{i,\sigma(j)} } \frac{{ b_j }}{{ a_{j,\sigma(j)} } }.\quad (1) $$

Now, note that if $j\not=j'$,   we should have $\sigma(j)\not=\sigma(j')$. Thus, $\sigma$ is a permutation defined on $[n]$. Hence, the first implication is established. To prove the converse note that, condition $(1)$ implies that for all $j$
$$x^\star_{\sigma(j)}=\frac{b_j}{a_{j,\sigma(j)}}.$$
We have for all $i$,  for all $\max_{k\in [n]}\big(a_{i,\sigma(k)}x_{\sigma(k)}^\star\big)=\max_{k\in [n]}\big(a_{i,\sigma(k)}\frac{b_k}{a_{k,\sigma(k)}}\big)=b_i$.
Consequently, $x^\star$ is a solution.
To end the proof, the strict inequality
$$b_i > \max_{j\in [n]\backslash \{i\}}{ a_{i,\sigma(j)} } \frac{{ b_j }}{{ a_{j,\sigma(j)} } }  \quad (2)$$
is equivalent to
$$\max_{j\in [n]}\big(a_{i,\sigma(j)}  x^\star_{\sigma(j)}\big)=a_{i,\sigma(i)}x^\star_{\sigma(i)}=b_i>\max_{j\in [n]\backslash \{i\}} \big(a_{i,\sigma(j)}x^\star_{\sigma(j)}\big)$$
for all $i\in [n]$. However, this latter condition is not compatible with the existence of some $u\leq x^\star$ such that $u_k<x_k^\star$ for some $k$, which ends the proof. $\Box $ \\

The next statement shows that it the limit system \eqref{IS} has a regular solution, then there exists a nonnegative solution to the system of maximum equations \eqref{MaxEq}.

\begin{lem}\label{unicregular}Let $A\in \mathcal M_n(\Real_+)$ be a square matrix such that $a_{i,j}\geq 0$ for all $i,j\in [n]$.  Suppose moreover that $b\in \Real_{+}^n$. Any solution of the limit system \eqref{IS} in $\Real_+^n$ is solution of system of maximum equations \eqref{MaxEq}. Moreover, if  the limit system has a regular solution $x^\star\in \Real^n$ then  the system of maximum equations \eqref{MaxEq} has a solution in $\Real_+^n$ that is $\sum_{i\in [n]}|x_i^\star|e_i$.
\end{lem}
{\bf Proof:} First note that if $x^\star \in \Real_+^n$ is solution of the limit system, then we have
  $b_i=\max_{j\in [n]}{a_{i,j}x_j}=b_i=\langle a_{i},x^\star\rangle_\infty^+ = \langle a_i,x^\star\rangle_\infty^-$. Hence $x^\star$ is solution of system \eqref{MaxEq}. Suppose now that $x^\star$ is a regular solution of the semi-continuous regularized system \eqref{IS}.
This implies that for all $i$, $b_i=\langle a_{i},x^\star\rangle_\infty^+ = \langle a_i,x^\star\rangle_\infty^-.$

Let us prove that  $y^\star=\sum_{i\in [n]}|x_i^\star|e_i$ is solution of system \eqref{MaxEq}. Let $J_{\circ}=\{j: x_{j}^\star<0\}$.
If $y^\star$ is not solution of system \eqref{MaxEq} then, since $A\in \mathcal M_n(\Real_+)$, there is some $i\in [n]$ and some $j_\circ\in J_\circ$ such that  $\langle a_i,y^\star\rangle_\infty^-=a_{i,j_\circ}| x_{j_\circ}^\star|>b_i$. However, this implies that
   $\langle a_i,x^\star\rangle_\infty^-=a_{i,j_\circ}  x_{j_\circ}^\star <0 \leq b_i$,
  which is a contradiction. Consequently $y^\star$ is solution of system \eqref{IS}. $\Box$\\

In the following, a condition is given to ensure that the Cramer
formula expressed in this idempotent and non-associative algebraic structure yields a solution to a system of maximum equations. This is  a limit case of the condition proposed by Kaykobad \cite{kayko}   when $p\longrightarrow \infty$.

\begin{prop}\label{semisyst2}Let $A\in \mathcal M_n(\Real)$ be a square matrix such that $a_{i,j}\geq 0$ for all $i,j\in [n]$. Suppose that $b\in \Real_{++}^n$.
 If there exists a permutation $\sigma:
[n]\longrightarrow [n]$ such that for all $i$, we have $a_{i,\sigma(i)}>0$ and $$b_i > \max_{j\in
[n]\backslash \{i\}}{ a_{i,\sigma(j)} } \frac{{ b_j }}{{
a_{j,\sigma(j)} } }, $$ then $|A|_\infty\not=0$. Moreover, there exists a
  solution $x^\star =\sum_{i\in
[n]}\frac{|A^{(i)}|_\infty}{|A|_\infty}e_i\in \Real_{+}^n$ to  the
system of maximum equations \eqref{MaxEq}.

Conversely, suppose that $|A|_\infty\not=0$. If $x^\star=\sum_{i\in
[n]}\frac{|A^{(i)}|_\infty}{|A|_\infty}e_i$ is a uniqueness regular solution of the limit system   \eqref{IS}  then $x^\star$ is a nonnegative solution of system of maximum equations \eqref{MaxEq}.
\end{prop}
{\bf Proof:}  We first
establish that $|A|_{\infty}\not=0$. Let $B={\mathrm{diag}}(b)$ be the
diagonal matrix such that for all $i\in [n]$, $B_{i,i}=b_i$. Since
$b\in \Real_{++}^n$, $B$ is $\varphi_p$-invertible for all $p$.
Moreover, for all $p$, $|B^{-1}|_p=\big(\prod_{i\in
[n)}b_i\big)^{-1}$. Let $A'\in \mathcal M_n(\Real_+)$ such that:
$$A'=B^{-1}A.$$
Since $|A'|_p=|B^{-1}_p||A|_p$ for all $p\in \mathbb N$, we deduce
taking the limit that: $|A'|_\infty=|B^{-1}|_\infty
|A|_\infty=\big(\prod_{i\in [n)}b_i\big)^{-1} |A|_\infty$. Hence,
$|A|_\infty\not=0$ if and only if $|A'|_\infty\not=0$. Since, $$b_i > \max_{j\in
[n]\backslash \{i\}}{ a_{i,\sigma(j)} } \frac{{ b_j }}{{
a_{j,\sigma(j)} } }, $$ we deduce that for all $i$ and all $j\not=i$:
$$a'_{j,\sigma(j) } > a'_{i,\sigma(j) }.$$
In particular this implies that for all $j$, $a'_{j,\sigma(j)
}>0$. From the limit form of the Leibniz formula, we deduce that:
$$|A'|_\infty=\mathrm{sgn}(\sigma)\prod_{j\in [n]}a'_{j,\sigma(j) }\not=0.$$
Hence $|A|_\infty\not=0$. Let us consider the system
$A\stackrel{p}{\cdot}x=b.$ We have established that if$$b_i > \Big(\sum_{j\in [n]\backslash
\{i\}}{(a_{i,j})}^{2p+1}\frac{{(b_j)}^{2p+1}}{{(a_{j,j})}^{2p+1}}\Big)^{\frac{1}{2p+1}}
$$for all
$i\in [n]$, then $A$ is $\varphi_p$-invertible and there is a solution
$x^{(p)}\in \Real_{++}^n$ to the system $A\stackrel{p}{\cdot}x=b$.
However, we have $$\lim_{p\longrightarrow \infty}\Big(\sum_{j\in
[n]\backslash \{i\}}{(a_{i,\sigma
(j)})}^{2p+1}\frac{{(b_j)}^{2p+1}}{{(a_{j,\sigma(j)})}^{2p+1}}\Big)^{\frac{1}{2p+1}}=\max_{j\in
[n]\backslash
\{i\}}\{a_{i,\sigma(j)}\frac{b_j}{a_{j,\sigma(j)}}\}.$$ Hence, there
is some $p_0\in \mathbb N$ such that for all $p\geq p_0$, we have
$|A|_p\not=0$ which implies that $x^{(p)}=\sum_{i\in
[n]}\frac{|A^{(i)}|_p}{|A|_p}e_i\in \Real_{++}^n$ is solution of the
system $A\stackrel{p}{\cdot}x=b$. However $x^\star=\sum_{i\in
[n]}\frac{|A^{(i)}|_\infty}{|A|_\infty}e_i=\lim_{p\longrightarrow
\infty}x^{(p)}$. It follows that $x^\star\in \Real_+^n$. We only
need to prove that for all $i\in [n]$,
$\max_{j}a_{i,j}x^\star_j=b_i.$ We have shown that $x^\star\in
Li_{p\longrightarrow \infty} \{x\in \Real^n:
A\stackrel{p}{\cdot}x=b\}$. Since $Li_{p\longrightarrow \infty}
\{x\in \Real^n: A\stackrel{p}{\cdot}x=b\}\subset \bigcap_{i\in
[n]}\big(\;[\langle a_i,\cdot\rangle_\infty^-\leq b_i]\cap
[\langle a_i,\cdot\rangle_\infty^+\geq b_i]\;\big)$. However since
$a_i\geq 0$ for all $i$, it follows that for all $i$:
\begin{equation}\langle a_i,x^\star\rangle_\infty^-=\langle a_i,x^\star\rangle_\infty^+=\max_ja_{i,j}x_j=b_i.\end{equation} Therefore
$x^\star$ is solution of   system of maximum equations. Conversely, if $|A|_\infty\not=0$, then $x^\star\in
Li_{p\longrightarrow \infty} \{x\in \Real^n:
A\stackrel{p}{\cdot}x=b\}$. Consequently, if $x^\star$ is regular, we deduce from  Lemma \ref{unicregular} that  $x^\star$ is a nonnegative solution  system of maximum equations
 \ref{MaxEq}. $\Box$\\

From Lemma \ref{eqforms} and Proposition \ref{semisyst2}, the following corollary is immediate.
\begin{cor}\label{eqforms2}Let $A\in \mathcal M_n(\Real)$ be a square matrix such that $a_{i,j}\geq 0$ for all $i,j\in [n]$.
For all $i,j\in [n]$, let us denote $\mu(i) = \{j : a_{i, j} > 0
\}$ and $\eta(j) = \{i : a_{i, j} > 0
\}$ and  assume that  $  \eta(j)$ and $ \mu(i)$ are nonempty. Suppose moreover that $b\in \Real_{++}^n$. If the system of maximum equations \eqref{MaxEq} has a   uniqueness solution in $\Real_+^n$    then $|A|_\infty\not=0$ and this solution is $x^\star =\sum_{i\in
[n]}\frac{|A^{(i)}|_\infty}{|A|_\infty}e_i\in \Real_{+}^n$.
\end{cor}

\begin{center}
 {\scriptsize
 \unitlength 0.4mm 
\linethickness{0.4pt}
\ifx\plotpoint\undefined\newsavebox{\plotpoint}\fi 
\begin{picture}(192.25,172.25)(0,0)
\put(87.25,16.75){\vector(0,1){147}}
\put(8.75,88.75){\vector(1,0){177.5}}
\put(169.75,132.75){\line(0,-1){.25}}
\put(192.25,89.5){\makebox(0,0)[cc]{$x_1$}}
\put(87,172.25){\makebox(0,0)[cc]{$x_2$}}
\put(79.5,81.75){\makebox(0,0)[cc]{$0$}}
\put(137.25,134.75){\makebox(0,0)[]{$x^\star$}}
\multiput(33.25,153.25)(.0336990596,-.0423197492){319}{\line(0,-1){.0423197492}}
\put(44,139.75){\line(1,0){85.5}}
\put(129.5,139.75){\line(0,-1){99.25}}
\multiput(129.5,40.5)(.0337338262,-.042051756){541}{\line(0,-1){.042051756}}
\put(85.5,0){\makebox(0,0)[cc]
{{\bf
Figure \ref{SystEq}.4} Positive Solutions of Limit Systems}}
\put(129.5,129.75){\circle*{1.5}}
\multiput(0,143.5)(.0612745098,-.0337009804){408}{\line(1,0){.0612745098}}
\put(25,129.75){\line(1,0){126.25}}
\put(151.25,129.75){\line(0,-1){79}}
\multiput(151.25,50.75)(.0490367776,-.0337127846){571}{\line(1,0){.0490367776}}
\end{picture}
}\end{center}
We illustrate these results on simple numerical examples.
\begin{expl} Let us consider the following system:
\begin{equation}\label{expS} \left\{
\begin{matrix}\max\{2 x_{1},3x_{2}\}&=1\\
\max\{4 x_{1},x_{2}\}&=1.\end{matrix} \right. \end{equation} We
have ${a_{1}}^{}=(2,3)^{}$, ${a_{2}}^{}=(4,1)^{}$, $b_1 = 1$ and
$b_2 = 1$, from which we get $u_1 = \min\left\{\frac{1}{2},
\frac{1}{4}\right\} = \frac{1}{4}\,  \hbox{ and } \, u_2 =
\min\left\{\frac{1}{3}, \frac{1}{1}\right\} = \frac{1}{3}$. One
can check that $(\frac{1}{4}, \frac{1}{3})$ is a solution of
system \eqref{expS}. Let us consider the matrices:
$$A=\begin{pmatrix}2&3\\4&1\end{pmatrix}\quad A^{(1)}=\begin{pmatrix}1&3\\1&1\end{pmatrix}\quad
A^{(2)}=\begin{pmatrix}2&1\\4&1\end{pmatrix}.$$
 We have $\left|A\right|_\infty=(2\cdot 1\boxplus (- 3\cdot
 4))=-12$; $\left|A^{(1)}\right|_\infty=(1\cdot 1\boxplus (- 3\cdot
1))=-3$; $\left|A^{(2)}\right|_\infty=(2\cdot
1\boxplus (- 1\cdot
 4))=-4$. One can then retrieve the above solutions:
 $$x_1^\star= \frac{1}{4}\quad x_2^\star
 =\frac{-4}{-12}=\frac{1}{3}.$$

\end{expl}

In the following a three dimensional example is given.

 \begin{expl}Let us consider the following system:
\begin{equation}\label{EqS2}\left\{
\begin{matrix}\max\{ x_{1},3x_{2}, 4 x_3\}&=1\\
\max\{2 x_{1},5x_{2}, x_3\}&=1\\ \max\{4 x_{1},2x_{2},
x_3\}&=1.\end{matrix} \right. \end{equation} We have
${a_{1}}^{}=(1,3,4)^{}$, ${a_{2}}^{}=(2,5,1)^{}$, $a_3=(4,2,1)$,
$b_1 = b_2 =b_3= 1$, from which we get $u_1 =
\min\left\{1,\frac{1}{2}, \frac{1}{4}\right\} = \frac{1}{4}$, $
u_2 = \min\left\{\frac{1}{2}, \frac{1}{5}, \frac{1}{3}\right\} =
\frac{1}{5}$, and $ u_3 = \min\left\{\frac{1}{4}, 1, 1\right\} =
\frac{1}{4}$. One can check that $(\frac{1}{4}, \frac{1}{5},
\frac{1}{4})$ is a solution of system \eqref{EqS2}. Let us
consider the matrices:
$$A=\begin{pmatrix}1&3&4\\2&5&1\\ 4&2&1\end{pmatrix}\quad A^{(1)}=\begin{pmatrix}1&3&4\\1&5&1\\ 1&2&1\end{pmatrix}\quad
A^{(2)}=\begin{pmatrix}1&1&4\\2&1&1\\ 4&1&1\end{pmatrix}\quad
A^{(3)}=\begin{pmatrix}1&3&1\\2&5&1\\ 4&2&1\end{pmatrix}. $$
 We have:\\
  $\left|A\right|_\infty=1\cdot 5\cdot 1\boxplus  3\cdot
 1\cdot 4\boxplus 2\cdot 2\cdot 4 \boxplus (-4\cdot 5\cdot 4)\boxplus (-1\cdot 2\cdot 1)\boxplus
 (-2\cdot 3\cdot 1)=-80$; \\

\medskip
 \noindent  $\left|A^{(1)}\right|_\infty=1\cdot 5\cdot 1\boxplus  3\cdot
 1\cdot 1\boxplus 1\cdot 2\cdot 4 \boxplus (-4\cdot 5\cdot 1)\boxplus (-3\cdot 1\cdot 1)\boxplus (-2\cdot 1\cdot
 1)=-20$;\\

\medskip
\noindent  $\left|A^{(2)}\right|_\infty=1\cdot 1\cdot 1\boxplus  2\cdot
 1\cdot 4\boxplus 1\cdot 1\cdot 4 \boxplus (-4\cdot 1\cdot 4)\boxplus (-1\cdot 2\cdot 1)\boxplus
 (-1\cdot 1\cdot 1)=-16$ and
$\left|A^{(3)}\right|_\infty=1\cdot 5\cdot 1\boxplus  2\cdot
 2\cdot 1\boxplus 3\cdot 1\cdot 4 \boxplus (-4\cdot 5\cdot 1)\boxplus (-1\cdot 2\cdot 1)\boxplus
 (-2\cdot 3\cdot 1)=-20$. \\ One can then retrieve the above solutions:
 $$x_1^\star=\frac{-20}{-80}=\frac{1}{4}\quad x_2^\star =\frac{-16}{-80}=\frac{1}{5}\quad x_3^\star=\frac{-20}{-80}=\frac{1}{4}.$$

\end{expl}

\subsection{Limit Two-Sided Systems}

Let $A,C\in \mathcal M_n(\Real )$ and let $b,d\in \Real ^n$. We consider the following system:

\begin{equation}\label{twosided1}\left\{
 \begin{matrix}
 \big( A\stackrel{-}{\boxtimes}x\big)\stackrel{-}{\smile}d\leq \big(C\stackrel{-}{\boxtimes}x\big)\stackrel{-}{\smile}b&\\
   \big(A\stackrel{+}{\boxtimes}x\big)\stackrel{+}{\smile}d\geq \big(C\stackrel{+}{\boxtimes}x\big)\stackrel{+}{\smile}b,&x\in \Real^n\\
\end{matrix}\right. \end{equation}

In the following, we provide a sufficient condition for the existence of a solution and given.
To do that we introduce the matrix:
$$A\boxminus C=(a_{i,j}\boxminus c_{i,j})_{i,j\in [n]}$$
where the symbol $\boxminus$ means that for all $\alpha,\beta\in \Real$, $\alpha\boxminus \beta=\alpha\boxplus (-\beta)$.

\begin{prop}\label{semisyst2}Let $A,C\in \mathcal M_n(\Real )$ and let $b,d\in \Real ^n$. If $|A\boxminus B|_\infty\not=0$, then
$$x^\star=\sum_{i\in [n]} \frac{|(A\boxminus B)^{(i)}|_\infty}{| A\boxminus B  |_\infty}e_i$$
is solution of system \eqref{twosided1}, where $(A\boxminus B)^{(i)}$ is the matrix obtained by replacing the $i$-th column with $b\boxminus d$. Moreover, $\{x^\star\}= Lim_{p\longrightarrow \infty}\big \{x\in \Real^n :(A\boxminus B)\stackrel{p}{\cdot}x=(b\boxminus d) \big \}. $ It follows that $x^\star$ is a solution of the limit system:
\begin{equation}\Bigg\{\begin{matrix}& \big(A\boxminus C\big) \stackrel{-}{\boxtimes} x&\leq& b\boxminus d&\\& \big(A\boxminus C\big) \stackrel{+}{\boxtimes} x&\geq& b\boxminus d,& x\in \Real^n.\end{matrix}\end{equation}

\end{prop}

\noindent {\bf Proof:} Let $(a\boxminus c)_i=a_i\boxminus c_i$ denotes the $i$-th line of the matrix $A\boxminus C$. Moreover, for all natural numbers $p$, let us denote $a_i\stackrel{p}{-}c_i$ the $i$-th line of matrix $A\stackrel{p}{-}C$.  We have $ a_i\boxminus c_i=\lim_{p\longrightarrow \infty}a_i\stackrel{p}{-}c_i$ and $b_i\boxminus d_i=\lim_{p\longrightarrow \infty}b_i\stackrel{p}{-}d_i$. This implies  from Proposition \ref{limhyper} that for all $i$:
 $$Li_{p\longrightarrow \infty}\big [\langle a_i\stackrel{p}{-}c_i,\cdot\rangle_p\leq b_i\stackrel{p}{-}d_i \big ]= Li_{p\longrightarrow \infty} \big [\langle  a_i\boxminus c_i,\cdot\rangle_p\leq b_i\boxminus d_i \big ]. $$
Moreover, we have:
$$ Li_{p\longrightarrow \infty}\bigcap_{i\in [n]}\big [\langle  a_i\boxminus c_i,\cdot\rangle_p= b_i\boxminus d_i \big ]\subset \bigcap_{i\in [n]} Li_{p\longrightarrow \infty}\big [\langle  a_i\boxminus c_i,\cdot\rangle_p= b_i\boxminus d_i \big ]. $$
Hence, we deduce that
$$Li_{p\longrightarrow \infty}\bigcap_{i\in [n]}\big [\langle  a_i\boxminus c_i,\cdot\rangle_p= b_i\boxminus d_i \big ]\subset \bigcap_{i\in [n]} Li_{p\longrightarrow \infty}\big [\langle a_i\stackrel{p}{-}c_i,\cdot\rangle_p\leq b_i\stackrel{p}{-}d_i \big ]. $$
Moreover, since $| A\boxminus B  |_\infty\not=0$, from Proposition \ref{semisyst}
$$x^\star=\sum_{i\in [n]} \frac{|(A\boxminus B)^{(i)}|_\infty}{| A\boxminus B  |_\infty}e_i \in Li_{p\longrightarrow \infty} \bigcap_{i\in [n]}\big [\langle  a_i\boxminus c_i,\cdot\rangle_p= b_i\boxminus d_i \big ].$$
Hence, we deduce that
\begin{equation}\label{inclus}x^\star  \in \bigcap_{i\in [n]} Li_{p\longrightarrow \infty}\big [\langle  a_i\stackrel{p}{-}c_i ,\cdot\rangle_p= b_i\stackrel{p}{-}d_i  \big ].\end{equation}
For all natural numbers $p$, let us denote:
$E_{i}^{(p)}= \{z\in \Real^{ n}\times \Real^{ n}\times \Real^2: \langle (a_i,-c_i, d_i,-b_i), z\rangle _p\leq 0]\},$
$F_1= \{z\in \Real^{n}\times \Real^{n}\times \Real^2: z_i=z_{i+n}:i\in [n]\}$ and $F_2= \{z\in \Real^{n}\times \Real^{n}\times \Real^2: z_{2n+1}=z_{2n+2}=1\}$.
However,
$$ \big\{(x,x,1,1)\in \Real^{2n+2}: \langle a_i\stackrel{p}{-}c_i,x\rangle_p\leq  b_i\stackrel{p}{-}d_i\big \}=E_i^{(p)}\cap F_1\cap F_2.
$$
Therefore
$$Li_{p\longrightarrow \infty} \big\{(x,x,1,1)\in \Real^{2n+2}: \langle a_i\stackrel{p}{-}c_i,x\rangle_p\leq  b_i\stackrel{p}{-}d_i\big \}\subset Li_{p\longrightarrow \infty} \big(E_i^{(p)}\cap  B_1\cap B_2\big).
$$
It follows that
$$z^\star=(x^\star,x^\star,1,1)\in \big(Li_{p\longrightarrow \infty}E_{ i }^{(p)}\big)\cap \big(F_1\cap F_2\big).$$
However
$$ Li_{p\longrightarrow \infty} E_{i}^{(p)}=\big [\langle (a_i,-c_i, d_i,-b_i), \cdot \rangle _\infty^-\leq 0\big]\cap
\big[\langle (a_i,-c_i, d_i,-b_i), \cdot \rangle _\infty^+\geq 0\big].$$
Hence:
\begin{align*}  \big(Li_{p\longrightarrow \infty} E_i^{(p)}\big)\cap \big(F_1\cap F_2\big) &= \\
\Big\{(x,x,1,1)\in \Real^{2n+2}:& \big(\bigsmileminus_{j\in [n]} a_{i,j}x_j \big) \stackrel{-}{\smile} \big( \bigsmileminus_{j\in [n]}(-c_{i,j})x_j \big) \stackrel{-}{\smile}   d_{i} \stackrel{-}{\smile}  ( -b_{i})\leq 0, \\&\big(\bigsmileplus_{j\in [n]}  a_{i,j} x_{j}\big) \stackrel{+}{\smile}\big( \bigsmileplus_{j\in [n]}(-c_{i,j})x_j \big) \stackrel{+}{\smile}   d_{i} \stackrel{+}{\smile}  ( -b_{i})\geq 0\Big \}\end{align*}
Now, note that, for all real numbers  $\alpha,\beta$ $$\alpha\leq \beta\iff \alpha\stackrel{-}{\smile}(-\beta)\leq 0 \iff 0 \leq (-\alpha)\stackrel{+}{\smile}\beta.$$
Starting from System \eqref{twosided1}, we have for all $i\in [n]$:
\begin{align*}  & \big(\bigsmileminus_{j\in [n]} a_{i,j}x_j \big) \stackrel{-}{\smile} \big( \bigsmileminus_{j\in [n]}(-c_{i,j})x_j \big) \stackrel{-}{\smile}   d_{i} \stackrel{-}{\smile}  ( -b_{i})\leq 0
\notag \\ & \iff   \big( \bigsmileminus_{j\in [n]}a_{i,j}x_j\big) \stackrel{-}{\smile}   d_{i} \leq  \big( \bigsmileminus_{j\in [n]} c_{i,j} x_j \big) \stackrel{-}{\smile}   b_{i} \end{align*}
and
\begin{align*} &\big(\bigsmileplus_{j\in [n]}  a_{i,j} x_{j}\big) \stackrel{+}{\smile}\big( \bigsmileplus_{j\in [n]}(-c_{i,j})x_j \big) \stackrel{+}{\smile}   d_{i} \stackrel{+}{\smile}  ( -b_{i})\geq 0 \\ &\iff \big(\bigsmileplus_{j\in [n]} a_{i,j}x_j \big)
\stackrel{+}{\smile}   d_{i} \geq \big( \bigsmileplus_{j\in [n]} c_{i,j} x_j \big) \stackrel{+}{\smile}   b_{i}. \notag \end{align*}
 Hence from equation \eqref{inclus}, and since $(x^\star,x^\star,1,1)\in \bigcap_{i\in [n]}\big(Li_{p\longrightarrow \infty} E_i^{(p)}\big)\cap \big(F_1\cap F_2\big)$ we deduce that $x^\star$ satisfies system \eqref{twosided1}. $\Box$\\

 If the   the matrices $A=(a_{i,j})_{i,j\in [n]}$, $C=(c_{i,j})_{i,j\in [n]}$ and the vectors $b$ and $d$ have positive entries, the problem of finding a nonnegative solution to the system \eqref{twosided1} can be written:
\begin{equation}\label{twosided}\left\{
\begin{matrix}
\max\{a_{1, 1}x_{1},\ldots, a_{1, n}x_{n}, d_1\}&=&\max\{c_{1, 1}x_{1},\ldots, c_{1, n}x_{n}, b_1\}\\
\vdots&\vdots&\vdots \\
\max\{a_{m, 1}x_{1},\ldots, a_{m, n}x_{n}, d_n\}&=&\max\{c_{m, 1}x_{1},\ldots, c_{m, n}x_{n}, b_m\}.\\
\end{matrix}\right.\end{equation}

A solution of system  \eqref{twosided1} is said to be {\bf regular} if for all $i\in [n]$:
\begin{equation}\bigsmileminus_{j\in [n]}a_{i,j}x_j^\star \stackrel{-}{\smile}d_i=\bigsmileplus_{j\in [n]}a_{i,j}x_j^\star \stackrel{+}{\smile}d_i\quad \text{and}\quad \bigsmileminus_{j\in [n]}c_{i,j}x_j^\star \stackrel{-}{\smile}b_i=\bigsmileplus_{j\in [n]}c_{i,j}x_j^\star \stackrel{+}{\smile}b_i.\end{equation}

 \begin{prop}\label{semisyst2}Let $A,C\in \mathcal M_n(\Real_+ )$ and let $b,d\in \Real_+ ^n$.
 If $x^\star$ is a regular solution of system \eqref{twosided1} then it is solution of system \eqref{twosided}, moreover
 $ \sum_{i\in [n]}|x_i^\star|e_i$ is a nonnegative solution of \eqref{twosided}.
\end{prop}
{\bf Proof:} Suppose that $x^\star$ is a regular solution system \eqref{twosided1}. Let us denote $y^\star=\sum_{i\in [n]}|x_i|e_i$. For any equations $(i)$, we consider four cases:\\
$(i)$ $\bigsmileminus_{j\in [n]}a_{i,j}x_j^\star \stackrel{-}{\smile}d_i=b_i$. In such a case, since $b_i\geq 0$, $\bigsmileminus_{j\in [n]}a_{i,j}y_j^\star \stackrel{-}{\smile}d_i=b_i$\\
$(ii)$ $ b_i=\bigsmileminus_{j\in [n]}c_{i,j}x_j^\star \stackrel{-}{\smile}b_i$. Similarly, since $d_i\geq 0$, $\bigsmileminus_{j\in [n]}c_{i,j} y_i^\star \stackrel{-}{\smile}b_i=d_i$\\
$(iii)$ Suppose that $(i)$ and $(ii)$ do not holds. In such a case:
$$ \bigsmileminus_{j\in [n]}a_{i,j}x_j^\star \stackrel{-}{\smile}d_i=\bigsmileminus_{j\in [n]}a_{i,j}x_j^\star =\bigsmileminus_{j\in [n]}c_{i,j}x_j^\star =\bigsmileminus_{j\in [n]}c_{i,j}x_j^\star \stackrel{-}{\smile}b_i.$$If $\bigsmileminus_{j\in [n]}a_{i,j}x_j^\star =\bigsmileminus_{j\in [n]}c_{i,j}x_j^\star<0$, then there is some $j_0,k_0\in [n]$ such that $x_{j_0}^\star<0$, $x_{k_0}^\star<0$ and such that
$$a_{i,j_0}x_{j_0}=\bigsmileminus_{ji\in [n]}a_{i,j}x_j^\star =\bigsmileminus_{j\in [n]}c_{i,j}x_j^\star=c_{i,k_0}x_{k_0}. $$
It follows that
$$\bigsmileminus_{j\in [n]}a_{i,j}y_j^\star =- a_{i,j_0}x_{j_0}^\star=- c_{i,k_0}x_{k_0}^\star=\bigsmileminus_{j\in [n]}c_{i,j}y_j^\star>0,$$
which implies that
$$\bigsmileminus_{j\in [n]}a_{i,j}y_j^\star \stackrel{-}{\smile}d_i  =\bigsmileminus_{j\in [n]}c_{i,j}y_j^\star\stackrel{-}{\smile}b_i.$$
Since these properties hold for all $i$, we deduce the result. $\Box$\\

\begin{expl}Let us consider the  system
\begin{equation}\left\{\begin{matrix}\max\{2x_1,x_2,3\}&=\max\{ x_1,x_2,4\}\\
\max\{ x_1,3x_2,2\}&=\max\{ 2 x_1,2 x_2,3\}.\end{matrix}\right.\end{equation}
We have
$A=\begin{pmatrix}2&1\\1&3\end{pmatrix}$, $C=\begin{pmatrix}1&1\\2&2\end{pmatrix}$, $b=\begin{pmatrix}4\\3\end{pmatrix}$,  and $d=\begin{pmatrix}3\\2\end{pmatrix}$
$$A\boxminus C=\begin{pmatrix}2&0\\-2&3\end{pmatrix}\quad \text{ and }\quad b\boxminus d=\begin{pmatrix}4\\3\end{pmatrix}.$$
It follows that:
$$(A\boxminus C)^{(1)}=\begin{pmatrix}4&0\\3&3\end{pmatrix}\quad \text{ and }\quad (A\boxminus C)^{(2)}=\begin{pmatrix}2&4\\-2&3\end{pmatrix}.$$

We obtain
$$x_1^\star= \frac{12}{6}=2\quad\text{ and }\quad x_2^\star=\frac{8}{6}=\frac{4}{3}.$$
\end{expl}

\subsection{Some Remarks on the Symmetrisation of Idempotent Semiring}\label{symetrization}The above algebraic structure can be viewed as some kind of non-associative symmetrization of the idempotent semi-ring $(\Real_+, \vee,\cdot)$. However there exist another approach  to construct a ring involving a balance relation and symmetrizing $(\Real_+, \vee, \cdot)$ (see \cite{h85} and
\cite{MP} in a Max-Plus context).
Following the usual construction of integers from natural numbers, one can introduce the following balance relation defined on $\Real_+^2\times \Real_+^2$ by:
\begin{equation}(x_+, x_-)\nabla (y_+, y_-)\iff  \max \{x_+, y_-\}= \max\{y_+, x_-\} , \end{equation}
where $x_+, x_-, y_+,y_-\in \Real_+ $. Let us denote $\mathbf{x}=(x_+,x_-)$ for all $(x_+,x_-)\in \Real_+^2$ and consider the quotient
$\mathbb S=\Real_+^2\backslash \nabla$. Let us define the operations
  $  \oplus  $ and  $ \otimes  $ on $\mathbb S$ as:
\begin{equation}
\mathbf{x}\oplus\mathbf{y}=(x_+\oplus y_+, x_-\oplus y_-)=(\max \{x_+, y_+\}, \max \{x_-, y_-\}),
\end{equation}
and
\begin{equation}
\mathbf{t}\otimes \mathbf{x}=(t_+  x_+\oplus t_-  x_-, t_+ x_-\oplus t_- x_+ ).
\end{equation}
$\mathbb S$ can be decomposed in three equivalence classes $\mathbb S_\oplus$, $\mathbb S_\ominus$ and $\mathbb S_\circ$ respectively associated to the sets  $ \{(x_+, t): t<x_+\}$ (called positive),
$\{(t, x_-): t<x_- \}$ (called negative) and  $\{(x_\circ , x_\circ)\}$ called balanced.  All the familiar identities valid in rings admit analogues replacing equalities
by balances. This means that associativity holds over $\mathbb S $. It follows that the binary operation $\oplus$ defined on $\mathbb S $ cannot be identified to the binary operation $\boxplus$. However, it can be related to the semi-continuous regularized operators $\stackrel{-}{\smile}$ and $\stackrel{+}{\smile}$.

Let $V: \mathbb S\longrightarrow \Real$ be the map defined as $V(\oplus x_+)=x_+$ for all $x_+\in \Real_+$, $V(\ominus x_-)=-x_-$ for all $x_-\in \Real_+$, and $V(x_\circ,x_\circ)=0$ for all $x_\circ\in \Real_+$. Suppose that
$(\mathbf{x}_1,...,\mathbf{x}_m)\in \mathbb S^m$. Then
\begin{equation}
V\Big(\bigoplus_{i\in [m]}\mathbf{x}_i\Big)=V\Big(\max_{i\in [n]}x_{i,+}, \max_{i\in [n]}x_{i,-}\Big)=\frac{1}{2}\Big(\bigsmileplus_{i\in [m]}V(\mathbf x_{i})+ \bigsmileminus_{i\in [m]}V(\mathbf x_{ i})\Big).
\end{equation}

Suppose that $\mathbf A=\big(\mathbf{a}_{i,j}\big)_{\substack{i=1...n\\j=1...n}}\in \mathcal M_n(\mathbb S)$. A determinant can be derived from this associative algebraic structure as:
\begin{equation}|\mathbf A|_{\mathbb S}=\bigoplus_{\sigma\in S_n}{\mathbf{sgn}}(\sigma) \bigotimes_{i\in [n]}{\mathbf{a}}_{i,\sigma(i)},\end{equation}
where $\mathbf{sgn}(\sigma)=\oplus 1$ if $\sigma$ is even and $\mathbf{sgn}(\sigma)=\ominus 1$ if $\sigma$ is odd.
Suppose that $A$ is a $3\times 3$-dimensional real matrix $$A=\begin{pmatrix}3&2&3\\1&3&2\\3&1&3\end{pmatrix}.$$
The positive components of $A$ can be identified to $\mathbb S_\oplus$. If $\mathbf{A}$ is the corresponding matrix, then
$|\mathbf A|_{\mathbb S}=(27,27)\nabla \mathbf{0}$ and we cannot derive a Cramer solution. However, one can check that
$|  A|_{\infty}=12\not=0$.

The symmetrization process described above is in general used in the context of Maslov's semi module where we replace $\vee$ with $\oplus$ and $+$ with $\otimes$ \cite{lms1}.   Applications can
be found in \cite{koloma} and \cite{ms} for
Max-Plus.  To be more precise, let $\mathbb M=\Real\cup \{-\infty\}$. For $x$ and $y$ in ${\mathbb M}^n$
let $d_{\hbox{\scriptsize M}_{+}}(x, y) = \mid\mid {\mathbf e}^x
-{\mathbf e}^y\mid\mid_{\infty}$ where ${\mathbf e}^x =
(e^{x_{1}}, \ldots, e^{x_{n}})$, with the convention $e^{-\infty}
= 0$, and, for $u\in{\mathbb R}_{+}^{n}$, $\mid\mid u \mid\mid =
\max_{1\in[ n]}x_i$. The map $x\mapsto {\mathbf e}^x$ is a
homeomorphism from ${\mathbb M}^n$ with the metric
$d_{\hbox{\scriptsize M}_{+}}$ to ${\mathbb R}_{+}^{n}$ endowed
with the metric induced by the norm $ \mid\mid \cdot
\mid\mid_{\infty}$; its inverse is the map $\mathbf{ln}(x) =
(\ln(x_1), \ldots, \ln(x_n) )$ from ${\mathbb R}_{+}^{n}$ to
${\mathbb M}^n$, with the convention $\ln(0) = -\infty$. For all $(t_1, \ldots,
t_n)\in [-\infty, 0]^n$ and all  $x_1, \ldots, x_n\in \mathbb M^n$, let us denote:
\begin{equation}
\bigoplus_{i=1}^nt_i\otimes x_i=\bigvee_{i=1}^{n}\big(x_i +t_i 1\!\!1_n\big).
\end{equation}
In the following a non-associative symmetrisation is proposed.  Suppose now
that $x\in \Real_-$ and let us extend the logarithm function to
the whole set of real numbers. This we do by introducing the set
\begin{equation}\widetilde {\mathbb M}=\mathbb M \cup (\Real +i\pi)\end{equation}
where $i$ is the complex number such that $i^2=-1$ and $\Real
+i\pi=\{x+i\pi:x\in \Real\}$. In the following we extend the
logarithmic function to $\widetilde {\mathbb M}$.  $\psi_{\ln}
:{\mathbb M} \longrightarrow \widetilde {\mathbb M}$ defined by
\begin{equation}
\psi_{ \ln}(x)=\left\{\begin{matrix}\ln(x)&\text{ if
}x>0\\-\infty&\text{ if }x=0\\\ln(-x)+i\pi&\text{ if
}x<0.\end{matrix}\right.
\end{equation}
The map $x\mapsto \psi_{ \ln}(x)$  is an isomorphism from $\mathbb
M$ to $\widetilde {\mathbb M}$. Let $\psi_{ \exp}(x): \widetilde
{\mathbb M}\longrightarrow \mathbb M$ be its inverse. Notice that
$\psi_{ \ln}(-1)=i\pi$. The scalar multiplication is is extended
to the binary operation $\widetilde{\otimes}:\widetilde{\mathbb
M}\times \widetilde{\mathbb M}\longrightarrow \widetilde{\mathbb
M}$ defined by
\begin{equation}
\left\{\begin{matrix}x&\widetilde{\otimes}&y&=&y
&\widetilde{\otimes}&x&=&x+y\\x&\widetilde{\otimes}&(
y+i\pi)&=&(y+i\pi)
&\widetilde{\otimes}&x&=&x+y+i\pi\\

(x+i\pi )&\widetilde{\otimes} &(y+i\pi)&=&(y+i\pi)
&\widetilde{\otimes}&(x+i\pi)&=&x+y\\

(x+i\pi )&\widetilde{\otimes}& -\infty&=&-\infty
&\widetilde{\otimes}&(x+i\pi)&=&-\infty.
\end{matrix}\right.
\end{equation}

For all $z\in \widetilde{\mathbb M}$ the symmetrical element is
$\tilde z=i\pi \otimes z$.
One can then introduce a corresponding absolute value function
$|\cdot |_{\widetilde {\mathbb M}}: \widetilde{\mathbb  M}\longrightarrow \mathbb R\cup \{-\infty\}$
defined by:

\begin{equation}
|z|_{\widetilde {\mathbb M}}=\left\{\begin{matrix}z-i\pi&\text{ if }& z\in \Real +i\pi\\z
&\text{ if }& z\in \Real\\-\infty&\text{ if }&
z=-\infty.\end{matrix}\right.
\end{equation}
This absolute value allows us to define the following binary
operation  on $\widetilde{\mathbb M}\times \widetilde{\mathbb M}$:
\begin{equation}z\widetilde{\boxplus} u=
\left\{\begin{matrix}z\ &\hbox{ if } &|z|_{\widetilde {\mathbb M}}&>&|u|_{\widetilde {\mathbb M}}\\
z&\hbox{ if }&z&=&u\\
-\infty &\hbox{ if }& \tilde z&=&u\\
u& \hbox{ if }& |z|_{\widetilde {\mathbb M}}&<&|u|_{\widetilde {\mathbb M}}.\end{matrix}\right.\end{equation} By
definition we have $z\widetilde{\boxplus}
u=\psi_{\ln}\big(\psi_{\exp}(z)\boxplus \psi_{\exp}(u)\big )$. Moreover, we have $z\widetilde{\otimes}
u=\psi_{\ln}\big(\psi_{\exp}(z)\otimes \psi_{\exp}(u)\big )$. For all $z\in \widetilde{M}^n$, let  us  denote \begin{equation}\widetilde{\bigboxplus}_{i\in [n]}z_i=
\psi_{\ln}\big(\bigboxplus_{i\in [n]}\psi_{\exp}(z_i)\big ).\end{equation}

In the remainder, we introduce an sign function $\widetilde{\sgn}$ defined on $S_p$ such that $\widetilde{\sgn}(\sigma)=1$ if $\sigma$ is even and  $\widetilde{\sgn}(\sigma)=i\pi$ if $\sigma$ is odd.  Suppose that $A$ is a square matrix of $\mathcal M_n(\widetilde {\mathbb M})$. The  symmetrized determinant defined on $\widetilde {\mathbb M}$ is now:
\begin{equation} |A|_{\widetilde {\mathbb M},\infty}=\psi_{\ln} \big(|\psi_{\exp}(A)|_\infty\big)=\bigboxplus_{\sigma\in S_n} \big(\widetilde {\sgn}(\sigma)
{\widetilde{ \medotimes}}_{i\in [n]}a_{i,\sigma(i)}\big).\end{equation}

\section{Eigenvalues in Limit}
In the following, we say that $\lambda\in \Real $ is {\bf an eigenvalue of $A$ in limit}, if:  $(1)$ there exists a sequence $\{(\lambda_p, v_p)\}_{p\in \mathbb N}\subset \Real\times \Real^n$  such that for all $p$, $A\stackrel{p}{\cdot}v_p=\lambda_pv_p$; (2) there is an increasing  subsequence $\{p_k\}_{k\in \mathbb N}$ with $\lim_{k\longrightarrow \infty}(\lambda_{p_k},v_{p_k})=(\lambda, v)$. $v$ is called  { \bf an eigenvector in limit of $A$}.

We start with the following intermediary result which will be
useful in the following. We say that for all
$\lambda \in \Real$, {\bf $P_A^{(p)}(\lambda)=|A-\lambda I|_p$ is
a $\varphi_p$-characteristic polynomial in $\lambda$}.

\begin{prop}Let $A\in \mathcal M_n(\Real)$ be a square matrix. Let $\lambda \in \Real$ and let us consider the matrix $A\stackrel{p}{-}\lambda I_n$ where
$I_n$ is the $n$-dimensional identity matrix. Then the
$\varphi_p$-polynomial $P_A^{(p)}(\lambda)$ in $\lambda$ is
$$P_A^{(p)}(\lambda)=\stackrel{\varphi_p}{\sum_{k\in \{0\}\cup [n] }}(-1)^{n-k}  \!\!\!\!\!\!\stackrel{\varphi_p}{\sum_{ 1 \leq h_1 <\cdots  <
 h_k \leq  n}}\;\stackrel{\varphi_p}{ \sum_{\sigma \in S_{ h_1,\cdots, h_k}
 }} \!\!\!\!\!\!   \big(\mathrm{sign}(\sigma) \!\!\!\!\!\prod_{i\in \{h_1,...,h_k\} }
     \!\!\!{a_{i, \sigma(i)}} \big)\lambda^{n-k},$$
where $S_{ h_1,\ldots, h_k}$ denotes the set of all the
permutations
 defined on $ \{ h_1,\ldots, h_k\}$. Moreover for all $\lambda \in \Real$
\begin{align*}P_A^{(\infty)}(\lambda)=\lim_{p\longrightarrow \infty }P_A^{(p)}(\lambda)= \!\!\!\!\!\!\!\bigboxplus_{\substack{ 1 \leq h_1 <\cdots  <
 h_k \leq  n\\\sigma \in S_{ h_1,\ldots, h_k}, \, {k\in \{0\}\cup [n]}}} \!\!\!\!\! (-1)^{n-k}  \big(\mathrm{sign}(\sigma) \!\!\!\!\!\prod_{i\in \{h_1,...,h_k\} }
     {a_{i, \sigma(i)}} \big)\lambda^{n-k}.\end{align*}

\end{prop}
{\bf Proof:}
The first part of the statement is derived using the usual procedure making the formal substitution
$+\mapsto \stackrel{p}{+}$. Let us denote $q_n=\sum_{k=0}^nk! C_{n}^k$. Let $\mathcal B_{q_n}$
be the canonical basis of $\Real^{q_n}$ and let
$\{B_{k}\}_{k=0,...,n}$ be a partition of $\mathcal B_{q_n}$ such
that for all $k$,  $B_k=\{e_{ k,h_1,...,h_k,\sigma}: h\in
\{h_1,...,h_k\}, \sigma \in S_{ \{h_1,...,h_k\}} \}$. Hence, we
have $\mathrm{Card} B_k=k!C_n^k$. Let  \begin{equation}\label{vpol}\gamma_{A}=\sum_{k\in \{0\}\cup [n]}{\sum_{
1 \leq h_1 <\cdots <
 h_k \leq  n}} \sum_{\sigma \in S_{ h_1,\ldots, h_k}
 }   \big(\mathrm{sign}(\sigma) \!\!\!\!\!\prod_{i\in \{h_1,...,h_k\}
 } {a_{i, \sigma(i)}}  \big)e_{k,h,\sigma}.\end{equation} Let us introduce the
 transformation $ \tau_A: \Real\longrightarrow \Real^{q_n}$
 defined by
 \begin{equation}\label{taupol}\tau_A (\lambda)=\sum_{k\in \{0\}\cup [n]}{\sum_{
1 \leq h_1 <\cdots <
 h_k \leq  n}} \sum_{\sigma \in S_{ h_1,\ldots, h_k}
 }  \lambda^{n-k} e_{k,h,\sigma}.\end{equation}

 An elementary calculus shows that, for all $\lambda\in \Real$
 \begin{equation*}
P_A^{(p)}(\lambda)=\langle \gamma_A, \tau_A(\lambda)\rangle_p.
 \end{equation*}
 For all $u\in \Real^{q_n}$, we have $\lim_{p\longrightarrow
\infty}\langle \gamma_A, u\rangle_p=\langle \gamma_A, u\rangle_\infty$. Hence, $P_A^{\infty}(\lambda)=\lim_{p\longrightarrow
\infty}\langle \gamma_A, \tau_A(\lambda)\rangle_p$. $\Box$\\

$P_A^\infty$ is called the limit characteristic polynomial. Let us
introduce now the lower and upper characteristic polynomial,
respectively defined by

\begin{align}P_{A,-}^{(\infty)}(\lambda)=  {\bigsmileminus_{\substack{ 1 \leq h_1 <\cdots  <
 h_k \leq  n\\\sigma \in S_{ h_1,\ldots, h_k}, \, {k\in \{0\}\cup [n]}}}}  (-1)^{n-k}  \big(\mathrm{sign}(\sigma) \!\!\!\!\!\prod_{i\in \{h_1,...,h_k\} }
     {a_{i, \sigma(i)}} \big)\lambda^{n-k}\end{align}
     and
     \begin{align}P_{A,+}^{(\infty)}(\lambda)=  {\bigsmileplus_{\substack{ 1 \leq h_1 <\cdots  <
 h_k \leq  n\\\sigma \in S_{ h_1,\ldots, h_k}, \, {k\in \{0\}\cup [n]}}}}  (-1)^{n-k}  \big(\mathrm{sign}(\sigma) \!\!\!\!\!\prod_{i\in \{h_1,...,h_k\} }
     {a_{i, \sigma(i)}} \big)\lambda^{n-k}.\end{align}

\begin{prop} Let $A\in \mathcal M_n(\Real)$ be a square matrix.
We have:
\begin{align*}Lim_{p\longrightarrow \infty}[P_A^{(p)}=0]=[P_{A,-}^{\infty}\leq 0]\cap [P_{A,+}^{\infty}\geq 0]= [P_{A,-}^{\infty}+ P_{A,+}^{\infty}=0].\end{align*}
Moreover,
$\lambda \in \Real$ is an eigenvalue in limit if and only if:
\begin{align*}\lambda \in [P_{A,-}^{\infty}+ P_{A,+}^{\infty}=0].\end{align*}
\end{prop}
{\bf Proof:} Let us denote $\gamma_A$ and $\tau_A(\lambda)$ respectively as in equation \eqref{vpol} and \eqref{taupol}. From Proposition \ref{limhyper}, we have $Lim_{p\longrightarrow \infty}[\langle \gamma_A,\cdot\rangle_p=0]= [\langle \gamma_A,\cdot\rangle_\infty^-\leq 0]\cap [\langle \gamma_A,\cdot\rangle_\infty^+\geq 0]$. Since that map $\tau_A$ is continuous, $Lim_{p\longrightarrow \infty}[P_A^{(p)}=0]=Lim_{p\longrightarrow \infty}\big [\langle \gamma_A,\tau_A(\cdot)\rangle_p=0\big ]=\big [\langle \gamma_A,\tau_A(\cdot)\rangle_\infty^-\leq 0\big]\cap \big [\langle \gamma_A,\tau_A(\cdot)\rangle_\infty^+\geq 0\big ]$. Hence $Lim_{p\longrightarrow \infty}[P_A^{(p)}=0]=[P_{A,-}^{\infty}\leq 0]\cap [P_{A,+}^{\infty}\leq 0]$. The last equality is an immediate consequence of Lemma \ref{equal}.  The second part of the statement is immediate since from Proposition  \ref{limhyper}, $\lambda \in [P_{A,-}^{\infty}\leq 0]\cap [P_{A,+}^{\infty}\leq 0]$ if and only if there is an increasing sequence $\{p_q\}_{q\in \mathbb N}$ and a sequence of real numbers $\{\lambda_{p_q}\}_{q\in \mathbb N}$ such that $\lim_{q\longrightarrow \infty}\lambda_{p_q}=\lambda$ and $[P_A^{(p_q)}(\lambda_{p_q})=0]$ for all $q$.  $\Box$\\

Given a square matrix with positive entries  $A\in \mathcal
M_n(\Real_{++})$, the Perron-Frobenius theorem states that there is
an eigenvalue called the spectral radius of $A$ and denoted
$\rho_A$ such that $\rho_A\geq |\lambda|$ for all eigenvalues of $A$,
where $|\cdot|$ denotes the module of $\lambda$. $\rho_A$ is
related to an eigenvector $v_A\in \Real_{++}^n$, with
$Av_A=\rho_Av_A$. $\lambda>0$ is an eigenvalue in the sense of the
matrix product $\boxtimes$ (a $\boxtimes$-eigenvalue) if there is a positive vector $v\in
\Real_{+}^n$ such that $A\boxtimes v=\lambda v$. We say that
$\lambda$ is a $\varphi_p$-eigenvalue of $A$ if
$A\stackrel{p}{\cdot}v=\lambda v$ for some vector $v\in
\Real^n\backslash \{0\}$. If $A\in \mathcal M_n(\Real_{++})$, then
$\Phi_p(A)\in \mathcal M_n(\Real_{++})$. Hence $\Phi_p(A)$ is
endowed with a spectral radius $\rho_{\Phi_p(A)}$, and a vector
$u_A^{(p)}\in \Real_{++}^n$ such that
$\Phi_p(A)u_A^{(p)}=\rho_{\Phi(A)}u_A^{(p)}$. It follows that
setting $v_A^{(p)}=\phi_p^{-1}(u_A^{(p)})$  and
$\rho_A^{(p)}=\varphi_p^{-1}(\rho_{\Phi_p(A)})$ that
\begin{equation}A\stackrel{p}{\cdot}v_A^{(p)}=\rho_A^{(p)}v_A^{(p)}.\end{equation} In
such a case, $\rho_A^{(p)}$ is called a $\varphi_p$-eigenvalue of
$A$.
Notice that in the case where $A\in \mathcal M_n(\Real_{++})$ there is only one $\boxtimes$-eigenvalue in $\Real_{++}$ (see for instance \cite{b10}).

\begin{prop}\label{Perlim}Let $A\in \mathcal M_n(\Real_{++})$ be a square matrix. For all $\lambda \in
\Real_{++}$ and all vectors $v\in \Real_{+}^n\backslash \{0\}$
such that $A\boxtimes v=\lambda v$, there is an increasing
subsequence $\{p_q\}_{q\in \mathbb N}$ such that
$\lambda=\lim_{q\longrightarrow \infty} \rho_A^{(p_q)}$ and
$v=\lim_{q\longrightarrow \infty} v_A^{(p_q)}$ where for all $p$,
$A\stackrel{p}{\cdot}v_A^{(p)}=\rho_A^{(p)}v_A^{(p)}.$
\end{prop}
{\bf Proof:}  We first prove that the sequence of the $\varphi_p$
Perron-Frobenius eigenvalues converges to a
$\boxtimes$-eigenvalue. For all $p$, there is an upper bound of
$\rho_A^{(p)}$. Moreover, $v_A^{(p)}$ can be chosen so as
$\|v_A^{(p)}\|=1$. Hence, there exists a compact subset $K$ of
$\Real_+^{n+1}$ which contains the sequence $\{(\rho_A^{(p)},
v_A^{(p)})\}_{p\in\mathbb N}$. Therefore, one can extract an increasing sequence
$\{p_q\}_{q\in \mathbb N}$ such that there is some
$(\rho_A^\infty, v_A^\infty) \in \Real_+^{n+1}$ with
$\lambda=\lim_{q\longrightarrow \infty} \rho_A^{(p_q)}$ and
$v=\lim_{q\longrightarrow \infty} v_A^{(p_q)}$. Since for all $q$,
$v_A^{(p_q)}\in \Real_{++}^n$ and since $A\in \mathcal
M_n(\Real_{+})$, it follows that:
$$\lim_{q\longrightarrow \infty}A\stackrel{p_q}{\cdot}v_A^{(p_q)}=A\boxtimes
v_A^\infty=\rho_A^{\infty}v_A^{\infty}.$$ Since there exists an
uniqueness $\boxtimes$-eigenvalue, we deduce that
$\lambda=\rho_A^{\infty}$ and that $v_A^{\infty}$ is a
$\boxtimes$-eigenvector. $\Box$\\

\begin{prop}Let $A\in \mathcal M_n(\Real_{++})$ be a square
matrix. If $\lambda \in \Real_{++}$ is a $\boxtimes$-eigenvalue
then $\lambda\in [P_{A,-}^{\infty}+ P_{A,+}^{\infty}=0].$ Moreover, if $\lambda$
is maximal in $[P_{A,-}^{\infty}+ P_{A,+}^{\infty}=0]$, then it is a $\boxtimes$-eigenvalue and $\lambda=\lim_{q\longrightarrow \infty} \rho_A^{(p_q)}$.

\end{prop}
{\bf Proof:} From Proposition \ref{Perlim}, we have $\lambda
=\rho_A^\infty$, we deduce that $\lambda\in Lim_{p\longrightarrow
\infty}[P_A^{(p)}=0]=[P_{A,-}^{\infty}\leq 0] \cap
[P_{A,+}^{\infty}\geq 0]=[P_{A,-}^{\infty}+ P_{A,+}^{\infty}=0]$. Conversely, suppose that $\lambda\in
[P_{A,-}^{\infty}\leq 0] \cap [P_{A,+}^{\infty}\geq 0]$ and assume that $\lambda$ is maximal in $[P_{A,-}^{\infty}+ P_{A,+}^{\infty}=0]$.  Then
there is a real sequence $\{\lambda^{(p)}\}_{p\in \mathbb N}$ with
$\lambda^{(p)}\in [P_A^{(p)}=0]$ for all natural numbers $p$ and such
that $\lim_{p\longrightarrow \infty}\lambda^{(p)}=\lambda$. Let
$\{\rho_A^{(p)}\}_{p\in \mathbb N}$ such that $\rho_A^{(p)}$ is a
$\varphi_{p}$-Perron-Frobenius eigenvalue for all $p$. From
Proposition \ref{Perlim}, there is a $\boxtimes$-eigenvalue
$\mu$ such that $\mu=\lim_{p\longrightarrow
}\rho^{(p)}_A$. It follows that $\mu \in [P_{A,-}^{\infty}\leq
0] \cap [P_{A,+}^{\infty}\geq 0]=[P_{A,-}^{\infty}+ P_{A,+}^{\infty}=0]$. Suppose that $\mu\not=\lambda$ and let us show a contradiction. Since $\lambda$ is maximal, this  implies that $\mu<\lambda$. However, this also implies that  there is some $p_0\in \mathbb N$ such that for all
$p>p_0$, $\lambda^{(p)}>\rho_A^{(p)}$,
which is a contradiction. Consequently, $\mu= \lambda$ and it follows that $\lambda$ is a $\boxtimes$-eigenvalue, which ends the proof. $\Box$\\

\begin{expl}Let us consider the matrix
$\begin{pmatrix}2&1\\1&2\end{pmatrix}$. Clearly $2$ is a
$\boxtimes$-eigenvalue and   $v=(1,1)$ is a $\boxtimes$-eigenvector, since $A\boxtimes v= 2v. $ The $\varphi_p$ Perron-Frobenius eigenvalue is
$\rho_A^{(p)}=(2^{2p+1}+1^{2p+1})^{\frac{1}{2p+1}}$ and we have
$\lim_{p\longrightarrow \infty}\rho_A^{(p)}=2$. We have
$$P^{(p)}_A(\lambda)=\big((\lambda^{2})^{2p+1}-(2\lambda)^{2p+1}-
(2\lambda)^{2p+1}+{4}^{2p+1}-1\big)^{\frac{1}{2p+1}}.$$

  Hence, taking the limit
yield:
$$P^{\infty}_A(\lambda)=(\lambda^{2})\boxplus
(-2\lambda)\boxplus (-2\lambda)\boxplus {4}\boxplus (-1).$$
Therefore
$$P^{\infty}_{A,-}(\lambda)=(\lambda^{2})\stackrel{-}{\smile}
(-2\lambda)\stackrel{-}{\smile} (-2\lambda)\stackrel{-}{\smile}
{4}\stackrel{-}{\smile} (-1)$$ and
$$P^{\infty}_{A,+}(\lambda)=(\lambda^{2})\stackrel{+}{\smile}
(-2\lambda)\stackrel{+}{\smile} (-2\lambda)\stackrel{+}{\smile}
{4}\stackrel{+}{\smile} (-1).$$ We have
$P^{\infty}_{A,-}(2)=-4\leq 0$ and $P^{\infty}_{A,+}(2)=4\geq 0.$

\end{expl}

\begin{expl}Let us consider the matrix
$\begin{pmatrix}1&2&1\\2&2&9\\1&1&3\end{pmatrix}$. Clearly $3$ is
a $\boxtimes$-eigenvalue  and $v=(2,3,1)$ is a $\boxtimes$
eigenvector, since
$A\boxtimes v=3v$.  We
have \begin{align*}P^{(p)}_A(\lambda)=\Big(-(\lambda^{3})^{2p+1}
&+\big[(2 \lambda^2)^{2p+1}+(1 \lambda^2)^{2p+1}+(3
\lambda^2)^{2p+1}\big]\\&-\big[(2\cdot 3\cdot
\lambda)^{2p+1}-(1\cdot 9\cdot \lambda)^{2p+1} -(1\cdot 2\cdot
\lambda)^{2p+1}\\&\quad \quad +(2\cdot 2\cdot
\lambda)^{2p+1}-(1\cdot 3\cdot \lambda)^{2p+1}+(1\cdot 1\cdot
\lambda)^{2p+1}\big]\\&+\big[(1\cdot 2\cdot 3)^{2p+1}+(2\cdot
1\cdot 1)^{2p+1}+(2\cdot 9\cdot 1)^{2p+1}\\&\quad \quad-(1\cdot
2\cdot 1)^{2p+1} -(2\cdot 2\cdot 3)^{2p+1}+(1\cdot 9\cdot
1)^{2p+1}\big]\Big)^{\frac{1}{2p+1}}.\end{align*}

  Hence, taking the limit
yield:

\begin{align*}P^{\infty}_A(\lambda)= -\lambda^{3}&\boxplus
2 \lambda^2\boxplus \lambda^2\boxplus 3 \lambda^2\\&\boxplus
6\lambda \boxplus (-9\lambda)\boxplus (- 2\lambda)\boxplus
4\lambda \boxplus (-3\lambda)\boxplus \lambda\\&\boxplus 3
\boxplus 2\boxplus 18\boxplus (-2)\boxplus  (-12)\boxplus
(-9).\end{align*}

Therefore

\begin{align*}P^{\infty}_{A,-}(\lambda)= -\lambda^{3} &\stackrel{-}{\smile}
2 \lambda^2\stackrel{-}{\smile} \lambda^2\stackrel{-}{\smile} 3
\lambda^2\\&\stackrel{-}{\smile} 6\lambda \stackrel{-}{\smile}
(-9\lambda)\stackrel{-}{\smile} (- 2\lambda)\stackrel{-}{\smile}
4\lambda \stackrel{-}{\smile} (-3\lambda)\stackrel{-}{\smile}
\lambda\\&\stackrel{-}{\smile} 3 \stackrel{-}{\smile}
2\stackrel{-}{\smile} 18\stackrel{-}{\smile}
(-2)\stackrel{-}{\smile} (-12)\stackrel{-}{\smile}
(-9).\end{align*} and
\begin{align*}P^{\infty}_{A,+}(\lambda)= -\lambda^{3}&\stackrel{+}{\smile}
2 \lambda^2\stackrel{+}{\smile} \lambda^2\stackrel{+}{\smile} 3
\lambda^2\\&\stackrel{+}{\smile} 6\lambda \stackrel{+}{\smile}
(-9\lambda)\stackrel{+}{\smile} (- 2\lambda)\stackrel{+}{\smile}
4\lambda \stackrel{+}{\smile} (-3\lambda)\stackrel{+}{\smile}
\lambda\\&\stackrel{+}{\smile} 3 \stackrel{+}{\smile}
2\stackrel{+}{\smile} 18\stackrel{+}{\smile} (-2)\boxplus
(-12)\stackrel{+}{\smile} (-9).\end{align*}We have
$P^{\infty}_{A,-}(3)=-27\leq 0$ and $P^{\infty}_{A,+}(3)=27\geq 0$

\end{expl}

In the next example, we consider a case where there is some $\mu
\in [P^{\infty}_{A,-}\leq 0] \cap [P^{\infty}_{A,+}\geq 0]=[P_{A,-}^{\infty}+ P_{A,+}^{\infty}=0]$ that
is an eigenvalue in limit but is not a $\boxtimes$-eigenvalue.
\begin{expl}Let us consider the matrix
$\begin{pmatrix}1&1\\1&1\end{pmatrix}$. $1$ is a
$\boxtimes$-eigenvalue associated to $v=(1,1)$ since
$A\boxtimes v=1.v$ The $\varphi_p$
Perron-Frobenius eignevalue is $\rho_A^{(p)}=2^{\frac{1}{2p+1}}$
and we have $\lim_{p\longrightarrow \infty}\rho_A^{(p)}=1$. For
all $p$, there is another eigenvalue $\mu^{(p)}=0$. We have
$$P^{(p)}_A(\lambda)=(\lambda^{2})^{2p+1}-\lambda^{2p+1}-
\lambda^{2p+1}.$$ Taking the limit yield:
$$P^{\infty}_A(\lambda)=(\lambda^{2})\boxplus
(-\lambda)\boxplus (-\lambda).$$ Therefore
$$P^{\infty}_{A,-}(\lambda)=(\lambda^{2})\stackrel{-}{\smile}
(-\lambda)\stackrel{-}{\smile} (-\lambda)\quad \text{ and }\quad P^{\infty}_{A,+}(\lambda)=(\lambda^{2})\stackrel{+}{\smile}
(-\lambda)\stackrel{+}{\smile} (-\lambda).$$
We have $P^{\infty}_{A,-}(1)=-1\leq 0$ and
$P^{\infty}_{A,+}(1)=1\geq 0$. There is another solution $\mu=0$,
we have $P^{\infty}_{A,-}(0)=P^{\infty}_{A,+}(0)= 0$. However, $0$
is not a $\boxtimes$-eigenvalue.

\end{expl}



\end{document}